\documentclass[10pt]{report}
\usepackage{amsmath, amsthm, amssymb}
\usepackage[dvips]{graphicx}
\theoremstyle{theorem}
\newtheorem{theorem}{Theorem}[section]
\newtheorem{proposition}{Proposition}[section]
\newtheorem{lemma}{Lemma}[section]
\newtheorem{corollary}{Corollary}[section]

\theoremstyle{definition}
\newtheorem{example}[theorem]{Example}
\newtheorem{definition}[theorem]{Definition}
\newtheorem{remark}[theorem]{Remark}

\def \Proof{\noindent{\textbf{Proof:}$\quad$}}
\def \pdt{\frac{\partial}{\partial t}}
\def \pdtu{\frac{\partial u}{\partial t}}
\def \pdtv{\frac{\partial v}{\partial t}}
\def \pdtw{\frac{\partial w}{\partial t}}
\def \pdtp{\frac{\partial p}{\partial t}}
\def \la{\langle}
\def \ra{\rangle}
\def \bf{\textbf}
\def \up{\textup}

\begin{document}
\title{Stochastic Completeness of Graphs}
\author{Rados\l aw Krzysztof Wojciechowski}
\maketitle

\begin{abstract}
We analyze the stochastic completeness of a heat kernel on graphs which is a function of three variables: a pair of vertices and a continuous time,
for infinite, locally finite, connected graphs.  For general graphs,
a sufficient condition for stochastic completeness is given in terms
of the maximum valence on spheres about a fixed vertex.  That this result
is optimal is shown by studying a particular family of trees.  We also
prove a lower bound on the bottom of the spectrum for the discrete
Laplacian and use this lower bound to show that in certain cases the
Laplacian has empty essential spectrum.
\end{abstract}

\tableofcontents

\chapter{Introduction}
\section{Introduction and Statement of Results}
The purpose of this thesis is to model a diffusion process on
infinite graphs which is analogous to the flow of heat on an open
Riemannian manifold.  In particular, we are interested in the
stochastic completeness of this process and a precise borderline for when the
stochastic completeness breaks down. Stochastic completeness can be formulated
in several equivalent ways: as a property of the heat kernel, as the
uniqueness of bounded solutions for the heat equation, or as the
non-existence of bounded, non-zero, $\lambda$-harmonic (or
$\lambda$-subharmonic) functions for a negative constant $\lambda$.
In studying this property, we have benefited tremendously from the survey
article of Grigor'yan \cite{Grig} which discusses, in great depth,
 stochastic completeness in the case of Riemannian
manifolds. For graphs, the starting point for our work is the paper
of Dodziuk and Mathai \cite{Do-Mat} where it is shown that any graph whose valence is uniformly
bounded above by a constant is stochastically complete.

In the first part, we give a construction of the heat
kernel on a general graph via an exhaustion argument.  This is analogous
 to the construction on open Riemannian manifolds and we
follow the presentation given in \cite{Do2}.  We also point out that
one can construct the heat kernel by utilizing the spectral theorem
but the two constructions result in the same kernel \cite{Do2}.  Next, we
introduce the notion of stochastic completeness and prove the
equivalence of the various formulations mentioned above.  This
material is adapted from \cite{Grig}.  We then turn our focus to a class of
trees which we call \emph{model} because their definition is
analogous to that of rotationally symmetric or model manifolds.  The
defining property of these trees is that they contain a vertex
$x_0$, which we call the \emph{root} for the model, such that the valence at every other
vertex depends only on the distance from $x_0$.  Let $m(r)$ denote
this common number where $r$ denotes the distance from $x_0$. The
main result of this section says that such trees will be
stochastically complete if and only if
$\sum_{r=0}^{\infty}\frac{1}{m(r)} = \infty$.  We note here the
similarity between this criterion and the one for the recurrence of
the Brownian motion on a complete, model surface \cite{Milnor,
Grig}.

We then consider general trees and prove that if a tree contains a
stochastically incomplete model subtree then it must be
stochastically incomplete.  We first prove this in the special case
when the branching of the general tree is growing rapidly in all
directions from the root, then in the general case where the
branching is growing rapidly in just some direction.  Also, we show that if a tree is contained in a stochastically
complete model tree then it must be stochastically complete. The proof of this fact
follows from the more general statement that, for any graph $G$, if
there exists a vertex $x_0$ such that the maximum valence of vertices on
spheres centered at $x_0$ is not growing too rapidly then $G$ must
be stochastically complete.

Next, we prove theorems analogous to a result of Cheeger and
Yau \cite{Ch-Y} which compare the heat kernel on a model tree to
the heat kernel on a general tree.  Let $T_n$ be a model tree with
root vertex $x_0$ where $n(r)$ is $m(r)-1$, that is, one less than the
common valence of vertices on the sphere of radius $r$ about $x_0$.  We denote
the heat kernel on $T_n$ by $\rho_t(x_0,x)$ and first show that, as
a function of $x$, $\rho_t(x_0,x)$ only depends on the distance from
$x_0$. That is, if we let $r(x) = d(x,x_0)$, where $d(x,x_0)$
denotes the distance between $x$ and $x_0$, then we can write
$\rho_t(r) = \rho_t(x_0,r(x))$. Let $T$ denote a general tree with
heat kernel $p_t(x_0,x)$. Then, if the branching on $T$ is growing
faster in all directions from $x_0$ then the branching on
$T_n$, we show that $p_t(x_0,x) \leq \rho_t(r(x))$.  In a similar
fashion, if $T \subseteq T_n$, then $\rho_t(r(x)) \leq p_t(x_0,x)$.

We finish this chapter by considering an operator related to the combinatorial
Laplacian that we study throughout the rest of the thesis.  This operator, referred to here
as the \emph{bounded Laplacian}, arises when one assigns the standard weight to the edges of a graph.  We show here that the
heat kernel associated to this Laplacian is stochastically complete for every graph $G$.  In particular,
bounded solutions for the combinatorial heat equation involving the bounded Laplacian are unique.

In the final part of the thesis we study the spectrum of the
Laplacian on a general graph. Specifically, we introduce
$\lambda_0(\Delta)$, the bottom of the spectrum of the Laplacian, and prove a
characterization of it in terms of the existence of positive
$\lambda$-harmonic functions.  That is, there always exist positive
functions satisfying $\Delta u = \lambda u$ for $\lambda \leq
\lambda_0(\Delta)$ whereas such functions never exist for $\lambda >
\lambda_0(\Delta)$ \cite{Sull, Do-Karp}.  We then prove a lower bound on
$\lambda_0(\Delta)$ under a geometric assumption on $G$.  Specifically, we
assume that if we fix a vertex $x_0$, then at every other vertex $x$ of $G$ the ratio of
the difference of the number of edges leaving $x$ and going away from $x_0$
and the number of edges going towards $x_0$ divided by the total valence
at $x$ is bounded below by a positive constant. The lower bound is then given in
terms of this constant.  In the final
section, we use this lower bound to prove that, with the additional
assumption that the minimum valence on spheres about $x_0$ is going to infinity
as one moves away from the fixed vertex, the Laplacian on the graph has empty
essential spectrum. This is analogous to the result of Donnelly and
Li for complete, simply connected, negatively curved Riemannian
manifolds \cite{Don-Li}.

\subsection{Acknowledgments}  I would like to thank Professor J\'ozef Dodziuk, my advisor, who has assisted me at every stage of my doctoral studies.  Without his wisdom and patience much of this work would have not been completed.  I would also like to thank Professors Isaac Chavel and Leon Karp for serving on my committee and for taking the time to discuss my results with me.


\section{Notation and Fundamentals}
In this section we fix our notation and state and prove some basic
lemmas which will be used throughout.  In general, $G=(V,E)$ will
denote an infinite, locally finite, connected graph where $V=V(G)$
is the set of vertices of $G$ and $E=E(G)$ the set of edges.  At
times, we abuse notation and write $x \in G$ when $x$ is a vertex of
$G$. We will use the notation $x \sim y$ to indicate that an edge
connects the vertices $x$ and $y$ while $[x,y]$ will denote the
\emph{oriented} edge from $x$ to $y$. In general, to be able to
write down certain formulas, we will assume that our graphs come
with an orientation, that is, that every edge is oriented, but none
of our results depend on the choice of this orientation. We use the
notation $m(x)$ to indicate the \emph{valence} at a vertex $x$, that
is, the number of edges emanating from $x$.

For a finite subgraph $D$ of $G$, we let Vol($D$) denote the \emph{volume}
of $D$ which we take, by definition, to be the number of vertices of
$D$.  That is,
\[ \textrm{Vol}(D) = \# \{x \ | \ x \in V(D) \}. \]
We also use the usual notion of distance between two vertices
of the graph.  Specifically, $d(x,y)$ will denote the number of edges in
the shortest path connecting the vertices $x$ and $y$.

We call $f$ a \emph{function on the graph} $G$ if it is a mapping
$f: V \to \mathbf{R}$.  The set of all such functions will be
denoted by $C(V)$. We will also use the notation $C_0(V)$ for the
space of all finitely supported functions on $G$ and $\ell^2(V)$ for
the space of all square summable functions.  That is, $\ell^2(V)$
consists of all functions on $G$ which satisfy
\[ \sum_{x \in V} f(x)^2 < \infty \]
and is a Hilbert space with inner product
\[ \la f, g \ra = \sum_{x \in V} f(x) g(x). \]
Similarly, we let $\ell^2(\tilde{E})$ denote the Hilbert space
of all square summable functions on oriented edges satisfying the
relation $\varphi([x,y]) = -\varphi([y,x])$ with inner product
\[ \la \varphi, \psi \ra = \sum_{[x,y] \in \tilde{E}} \varphi([x,y]) \psi([x,y]) \]
where $\tilde{E}$ denotes the set of all oriented edges of $G$.

We now recall the definitions of the \emph{coboundary} and
\emph{Laplacian} operators and state and prove an analogue of
Green's Theorem for them. The coboundary operator $d$ takes a
function on the vertices of $G$ and sends it to a function on the
oriented edges of $G$ defined by:
\[ df([x,y]) = f(y) - f(x) . \]
The combinatorial Laplacian $\Delta$ operates on functions on $G$ by the formula:
\begin{equation}\label{Laplacian}
\Delta f(x) = \sum_{y \sim x}\big(f(x) - f(y)\big) = m(x) f(x) -
\sum_{y \sim x} f(y)
\end{equation}
where the summation is taken over all vertices $y$ such that $y \sim
x$ forms an edge in $G$. If the Laplacian is applied to a function of
more than one variable then we will put the variable in which it is
applied as a subscript when necessary. For a constant $\lambda$, we
call a function $v$ on $G$ \emph{$\lambda$-harmonic} if
$\Delta v(x) = \lambda v(x)$ for all vertices $x$.

Note that it follows from formula (\ref{Laplacian}) that the
Laplacian will be bounded if and only if there exists a constant $M$ such
that $m(x) \leq M$ for all vertices $x$. Indeed, letting $\delta_x$
denote the delta function at a vertex $x$ so that
\[ \delta_x(y) = \left\{  \begin{array}{ll}
 1 & \textrm{ if } x = y\\
 0 & \textrm{ otherwise} \end{array} \right.
\]
we see that the matrix coefficients of the Laplacian are given
by
\begin{eqnarray*}
\Delta(x,y) &=& \la \Delta \delta_x, \delta_y \ra = \Delta \delta_x(y) \\
&=& \left \{ \begin{array}{lll} m(x) & \textrm{ if } x = y \\
-1 & \textrm{ if } x \sim y \\
0 & \textrm{ otherwise. } \end{array} \right.
\end{eqnarray*}
As mentioned in the introduction, under the assumption $m(x) \leq M$, all graphs are stochastically complete \cite[Theorem 2.10]{Do-Mat}.  Therefore, for the
purposes of our inquiry, we do not impose this restriction on the graph and the
Laplacian will be an unbounded operator.

Let $D$ be a finite, connected subgraph of $G$.  We then have
the following analogue of Green's Theorem.
\begin{lemma}\label{Green's}
\begin{eqnarray*}
\sum_{x \in D} \Delta f(x) g(x) &=& \sum_{[x,y] \in \tilde{E}(D)} df([x,y])
dg([x,y]) +  \sum_{\substack{x \in D \\ z \sim x, z \not \in D}}
\big(f(x) - f(z)\big)
g(x) \\
&=& \sum_{[x,y] \in \tilde{E}(D)} df([x,y]) dg([x,y]) -
\sum_{\substack{[x,z] \\ x \in D, z \not \in D}} df([x,z]) g(x).
\end{eqnarray*}
\end{lemma}
\Proof  Every oriented edge $[x,y]$ with $x,y \in V(D)$ contributes
two terms to the sum on the left hand side:
$\big(f(x)-f(y)\big)g(x)$ from $\Delta f(x)g(x)$ and
$\big(f(y)-f(x)\big)g(y)$ from $\Delta f(y)g(y)$. These add up to
give
\[ \big(f(y)-f(x)\big)\big(g(y)-g(x)\big) = df([x,y])dg([x,y]). \]
The remaining contributions come from any vertex $x$ in $D$ that is
connected to a neighbor $z$ which is not in $D$ and these give
$\big(f(x)-f(z)\big)g(x) = -df([x,z])g(x)$. \qed \\

We say that a vertex $x$ is in the \emph{boundary} of $D$, and
denote this $x \in \partial D$, if it is a vertex of $D$ and is
connected to any vertex which is not in $D$.  Otherwise, a vertex
$x$ of $D$ is said to be in the \emph{interior} of $D$, or $x \in$
int$(D)$.   We then see that, if either $f$ or $g$ are zero on the complement of
the interior of $D$, then the second term on the right hand side of the
equation above is zero and we can write Lemma \ref{Green's} as
\[ \la \Delta f, g \ra_{V(D)} = \la df, dg \ra_{\tilde{E}(D)} =
\la f, \Delta g \ra_{V(D)}. \]
Also, if $f$ and $g$ are any two functions and one of them is
finitely supported, it is true that
\[ \la \Delta f, g \ra = \la df, dg \ra = \la f, \Delta g \ra \]
where now the inner products are taken over $V(G)$ and
$\tilde{E}(G)$.

Throughout, we wish to study solutions of the \emph{combinatorial
heat equation}. These will be functions on $G$ with an additional
time parameter in which they are differentiable and which satisfy
the equation
\[ \Delta u(x,t) + \pdtu(x,t) = 0 \]
for every vertex $x$ and every $t > 0$.  We start by
recalling a proof of analogues for the weak and strong maximum principles for the
heat equation \cite{Jost, Do-Mat}.

\begin{lemma}\label{maxheatlem}
Suppose that $D$ is a finite, connected subgraph of $G$ and
\[u: D \times [0,T] \to \mathbf{R}\]
is continuous for $t \in [0,T]$, $C^1$ for $t \in (0,T)$, and
satisfies the combinatorial heat equation:
\[\Delta u + \pdtu = 0 \ on \ \textup{int } D \times (0,T).\]
Then, if there exists $(x_0,t_0) \in \textup{int} \ D \times (0,T)$
such that $(x_0,t_0)$ is a maximum (or minimum) for $u$ on $D \times
[0,T]$, then $u(x,t_0) = u(x_0,t_0)$ for all $x \in D$.
\end{lemma}
\Proof At either a maximum or minimum, $\pdtu(x_0,t_0) = 0$, giving
that
\[\Delta u(x_0,t_0) = \sum_{x \sim x_0}\Big(u(x_0,t_0) - u(x,t_0)
\Big) = 0.\] In either case, this implies that $u(x,t_0) =
u(x_0,t_0)$ for all  $x \sim x_0$. Iterating the argument and using the assumption that $D$ is connected gives the statement of the lemma. \qed

\begin{lemma}\label{maxheatcor}
Under the same hypotheses as above we have that
\[ \max_{D \times [0,T]} \ u = \max_{\substack{D \times \{0\} \ \cup\\
\partial D \times [0,T]}} \ u  \]
and
\[ \min_{D \times [0,T]} \ u = \min_{\substack{D \times \{0\} \ \cup \\
\partial D \times [0,T]}} \ u . \]
\end{lemma}

\Proof  Let $v = u - \epsilon t$ for $\epsilon > 0$.  Then $\Delta v + \pdtv = -\epsilon < 0.$  If $v$ has a maximum at $(x_0,t_0) \in \textrm{int} D \times (0,T]$ then
\[ \pdtv(x_0,t_0) \geq 0 \quad \textrm{and} \quad \Delta v(x_0,t_0) \geq 0 \]
yielding a contradiction.  Therefore,
\[ \max_{D \times [0,T]} \ v = \max_{\substack{D \times \{0\} \ \cup\\ \partial D \times [0,T]}} \ v. \]
Then
\begin{eqnarray*}
\max_{D \times [0,T]} \ u &=&  \max_{D \times [0,T]} \ v + \epsilon t  \\
&\leq&  \max_{D \times [0,T]} \ v + \epsilon T \\
&=&  \max_{\substack{D \times \{0\} \ \cup\\ \partial D \times [0,T]}} \ v + \epsilon T \\
&\leq&  \max_{\substack{D \times \{0\} \ \cup\\ \partial D \times [0,T]}} \ u + \epsilon T .
\end{eqnarray*}
Letting $\epsilon \to 0$ we get that
\[ \max_{D \times [0,T]} \ u = \max_{\substack{D \times \{0\} \ \cup\\
\partial D \times [0,T]}} \ u.  \]
The statement about the minimum follows by applying the argument to $-u$.    \qed

\begin{remark} \label{maxheatremark}
Using the same techniques as above, it follows that, if $u$
satisfies
\[\Delta u + \pdtu \geq 0 \ on \ \textup{int } D \times (0,T)\]
then
\[ \min_{D \times [0,T]} \ u = \min_{\substack{D \times \{0\} \ \cup\\
\partial D \times [0,T]}} \ u \]
while if $u$ satisfies
\[\Delta u + \pdtu \leq 0 \ on \ \textup{int } D \times (0,T) \]
then
\[ \max_{D \times [0,T]} \ u = \max_{\substack{D \times \{0\} \ \cup\\
\partial D \times [0,T]}} \ u. \] \\
\end{remark}

\section{Essential Self-Adjointness of the Laplacian}
As in the case of the Laplacian on a Riemannian manifold, the
Laplacian with domain $C_0(V)$, the set of all finitely supported
functions on the graph $G$, is a symmetric but not self-adjoint
operator. It is, however, essentially self-adjoint by which we mean
that it has a unique self-adjoint extension $\tilde{\Delta}$ to $\ell^2(V)$, a fact
which we prove in this section. Let $\Delta^*$ denote the adjoint of
$\Delta$ with domain $C_0(V)$.

\begin{proposition}
The domain of $\Delta^*$ is
\[ \textrm{dom}(\Delta^*) = \{ f \in \ell^2(V) \ | \ \Delta f \in \ell^2(V) \}. \]
\end{proposition}
\Proof  By definition
\[ \textrm{dom}(\Delta^*) = \left\{ \begin{array}{ll}
 f \in  \ell^2(V)  \ \left| \begin{array}{ll} & \textrm{there exists a unique } h \in \ell^2(V) \textrm{ such that } \\
 & \quad \la \Delta g, f \ra = \la g, h \ra \textrm{ for all } g \in
C_0(V) \end{array} \right.
\end{array} \right\}
\]
and then $\Delta^* f = h.$
If $g$ is finitely supported as above then we can apply
the analogue of Green's Theorem to get that if $f \in $
dom($\Delta^*)$ then
\[ \la \Delta g, f \ra = \la g, \Delta f \ra = \la g, h \ra. \]
Letting, $g = \delta_x$ we get that $\Delta f(x) = h(x)$ for all
vertices $x$ so that $\Delta f \in \ell^2(V)$. \qed

\begin{theorem}
$\Delta$ with domain $C_0(V)$ is essentially self-adjoint.
\end{theorem}
\Proof  From the criterion stated in \cite[Theorem
X.26]{Reed-SimonII} applied to the operator $(\Delta + I)$ it
suffices to show that $-1$ is not an eigenvalue of $\Delta^*$.  In
other words, if $f$ satisfies $\Delta^* f = - f$ then $f$ cannot be
in $\ell^2(V)$ unless it is exactly 0.  As can be seen by applying
the analogue of Green's Theorem, pointwise, it is true that
$\Delta^* f(x) = \Delta f(x)$.  Therefore, if $f$ satisfies
$\Delta^* f(x) = - f(x)$ for every vertex $x$, it follows that
\[ \big(m(x) + 1 \big) f(x) = \sum_{y \sim x} f(y). \]
Therefore, there must exist a neighbor $ y \sim x$ such that $f(y)
> f(x)$. By repeating this argument the conclusion follows.
\qed


\chapter{The Heat Kernel}

\section{Construction of the Heat Kernel}

We now give a construction of the heat kernel $p=p_t(x,y)$ for a infinite, locally finite,
connected graph $G$. By \emph{heat kernel} we mean that $p_t(x,y)$ will be
the smallest non-negative function
\[ p: V \times V \times \ [0,\infty) \to [0,1] \] which is smooth in $t$,
satisfies the heat equation: $\Delta p + \pdtp = 0$ in either $x$
or $y$ and satisfies: $p_0(x,y) = \delta_{x}(y)$.  The heat
kernel will generate a bounded solution of the heat equation on $G$ for any bounded initial condition.
That is, for any bounded function $u_0$,  $u(x,t) = \sum_{y \in V} p_t(x,y)u_0(y)$ will give a bounded
solution to
\[ \left\{  \begin{array}{ll}
 \Delta u(x,t) + \pdtu(x,t) = 0 & \textrm{ for all } x \in V, \textrm{all } t>0\\
 u(x,0) = u_0(x) & \textrm{ for all } x \in V.
 \end{array} \right. \]
The construction given here follows the approach of \cite[Section
3]{Do2} and the formalism of \cite{Do}.

Starting with an exhaustion sequence of the graph, we construct
heat kernels with Dirichlet boundary conditions for each set in the
exhaustion. Although we will exhaust the graph by balls of
increasing radii, it will be shown later that the resulting heat
kernel is independent of the choice of subgraphs in the exhaustion.
Let $x_0 \in V(G)$ be a fixed vertex. We will let $B_r = B_r(x_0)$
denote the ball of radius r about $x_0$, $\partial B_r =
\partial B_r(x_0)$ its boundary, and int $B_r$, its interior. In
particular, if $d(x,x_0)$ denotes the standard metric on graphs,
\begin{eqnarray*}
V(B_r) &=& \{ x \in V(G) \ | \ d(x,x_0) \leq r \} \\
E(B_r) &=& \{ x \sim y \ | \ x,y \in V(B_r) \textrm{ and } x \sim y
\in E(G) \}.
\end{eqnarray*}

We then let $C(B_r, \partial B_r)$ denote functions on $B_r$
which vanish on the boundary $\partial B_r$ and let $\Delta_r$
denote the \emph{reduced Laplacian} which acts on these spaces.
That is,
\[ C(B_r, \partial B_r) = \{f \in C(B_r) \ | \
f_{|\partial B_r} = 0 \} \] and
\[ \Delta_r f(x) = \left \{ \begin{array}{ll}
\Delta f(x) & \textrm{for } x \in \textrm{int } B_r\\
0 & \textrm{otherwise}
\end{array} \right.
\]
for all $f \in C(B_r, \partial B_r)$.

With these definitions we then have:
\begin{lemma}\label{positive}
$\Delta_r$ is a self-adjoint, non-negative operator on $C(B_r,
\partial B_r)$.
\end{lemma}
\Proof  This follows from the analogue of Green's Theorem, Lemma
\ref{Green's}, since
\[ \la \Delta_r f, g \ra_{V(B_r)} = \la df, dg \ra_{\tilde{E}(B_r)} =
\la f,\Delta_r g \ra_{B(B_r)}.
  \] \qed

From Lemma \ref{positive} it follows that all eigenvalues
$\lambda_i^r$ of $\Delta_r$ are real and non-negative. In fact, as
mentioned in \cite{Do}, and to be shown later,
$\lambda_0(\Delta_r) = \lambda_0^r$, the smallest eigenvalue of
$\Delta_r$, is given by the Rayleigh-Ritz quotient:
\[  \lambda_0^r = \min_{\substack{f \in C(B_r, \partial B_r) \\ f \not \equiv 0}}
\frac{\la df,df \ra}{\la f,f \ra} \] so that all of the eigenvalues
of $\Delta_r$ are positive.  Denote by $\{ \lambda_i^r
\}_{i=0}^{k(r)}$ the set of all eigenvalues of $\Delta_r$ listed in
increasing order and choose a set $\{ \phi_i^r\}_{i=0}^{k(r)}$ of
corresponding eigenfunctions which are an orthonormal basis for
$C(B_r, \partial B_r)$ with respect to the $\ell^2$ inner product.
That is, $\{ \phi_i^r \}_{0=1}^{k(r)}$ are such that,
\begin{equation} \label{eigen}
\Delta_r \phi_i^r = \lambda_i^r \phi_i^r \quad \forall i = 0, \dots,
k(r)
\end{equation}
and
\begin{equation} \label{on} \sum_{x \in B_r} \phi_i^r(x)
\phi_j^r(x) = \delta_{ij} =
\left \{ \begin{array}{ll} 1 & \textrm{if } i = j\\
0 & \textrm{otherwise.} \end{array} \right.
\end{equation}
We are now ready to define the heat kernels $p_t^r(x,y)$ for each
subgraph in the exhaustion.

\begin{definition}
\begin{equation}\label{ptr} p_t^r(x,y) = \sum_{i=0}^{k(r)} e^{-\lambda_i^r t}
\phi_i^r (x) \phi_i^r (y) \quad \textrm{ for all } x,y \in B_r,
\textrm{ all } t \geq 0.
\end{equation}
\end{definition}

\begin{theorem} \label{Dirichlet heat kernels}
$p_t^r(x,y)$ has the following properties for every $r$:
\begin{enumerate}
\item[\textup{1)}]  $p_t^r(x,y) = p_t^r(y,x)$, $p_t^r(x,y) = 0$ if
either $x \in \partial B_r$ or $y \in \partial B_r$.
\item[\textup{2)}]  $\Delta_r p_t^r(x,y) + \pdt p_t^r(x,y) = 0$
where $\Delta_r$ denotes the reduced Laplacian in either $x$ or $y$.
\item[\textup{3)}] $p_{s+t}^{r}(x,y) = \sum_{z \in B_r}p_s^r(x,z)
p_t^r(z,y)$.
\item[\textup{4)}] $p_0^r(x,y) = \delta_x(y)$ for $x,y \in \textup{int } B_r$.
\item[\textup{5)}] $p_t^r(x,y) > 0 \quad \textrm{ for all } t > 0,
\textrm{all } x,y \in \textup {int } B_r$.
\item[\textup{6)}] $\sum_{y \in B_r} p_t^r(x,y) < 1 \ \textrm{ for all } t >
0,  \textrm{all } x \in B_r$. \newline
\end{enumerate}
\end{theorem}

\begin{remark}
We could also start by defining the heat semigroup operator as the
convergent power series:
\[ Q_t^r = e^{-t\Delta_r} = I - t\Delta_r + \frac{t^2}{2}\Delta_r^2 -
\frac{t^3}{6}\Delta_r^3 + \dots \] and then take its kernel given by
$q_t^r(x,y) = \la Q_t^r \delta_x,\delta_y \ra = Q_t^r \delta_x(y)$.
The equivalence of these two approaches can be seen by applying the
maximum principle, Lemma \ref{maxheatcor}, to the difference of
the two kernels.
\end{remark}

\Proof \bf{1), 2)}  Clear from the definition of $p_t^r(x,y)$, the
fact that $\phi_i^r \in C(B_r, \partial B_r)$, that is, $\phi^r_{i
|\partial B_r} = 0$, and from (\ref{eigen}). \newline

\bf{3)}  Using the orthonormality of
$\{\phi_i^r\}_{i=1}^{k(r)}$, we compute:
\begin{eqnarray*}
\sum_{z \in B_r} p_s^r(x,z)p_t^r(z,y) &=& \sum_{z \in B_r}
\sum_{i=0}^{k(r)}e^{-\lambda_i^r s}\phi_i^r(x) \phi_i^r(z)
\sum_{j=0}^{k(r)}e^{-\lambda_j^r t} \phi_j^r(z) \phi_j^r(y)\\
&=& \sum_{i,j = 0}^{k(r)} e^{-\lambda_i^r s} e^{-\lambda_j^r t}
\phi_i^r(x) \phi_j^r(y) \sum_{z \in B_r} \phi_i^r(z) \phi_j^r(z) \\
&=& \sum_{i=0}^{k(r)} e^{-\lambda_i^r(s+t)} \phi_i^r(x) \phi_i^r(y)
\quad \quad \quad \quad \quad \quad \quad \quad \quad \quad (\ref{on})\\
&=& p_{s+t}^r(x,y).
\end{eqnarray*}

\bf{4)}  By definition,
\[ p_0^r(x,y) = \sum_{i=0}^{k(r)} \phi_i^r(x) \phi_i^r(y). \]
Since $\{ \phi_i^r \}_{i=0}^{k(r)}$ form an orthonormal basis, it follows that
\[ \delta_x(y) =  \sum_{i=0}^{k(r)} \la \delta_x , \phi_i^r \ra \phi_i^r(y) = \sum_{i=0}^{k(r)} \phi_i^r(x) \phi_i^r(y). \]
Therefore, $p_0(x,y) = \delta_x(y).$

\bf{5)}  The maximum principle, Lemma \ref{maxheatcor}, applied in each of the
variables separately to $p_t^r(x,y)$ over the set $B_r \times B_r \times [0,T]$
implies that $0 \leq p_t^r(x,y) \leq 1$ since $p_0^r(x,y) = \delta_x(y)$ and $p_t^r(x,y) =0$ if either $x$ or $y$ is in the boundary of $B_r$.

Now, assume that there exists a $t_0>0$ and $\hat{x},\hat{y}\in $
int $B_r$ such that $p_{t_0}^r(\hat{x},\hat{y}) = 0$.  We may then
assume that $(\hat{x},\hat{y},t_0)$ is a minimum for $p_t^r(x,y)$ on
$B_r \times B_r \times [0, t_0]$.  Then, using the fact that $p_t^r(x,y)$ satisfies the
heat equation in both variables and $B_r$ is connected, and applying the argument used in the proof of Lemma \ref{maxheatlem} over the set $B_r \times B_r \times [0, t_0]$ gives:
\[ p_{t_0}^r(x,y) = p_{t_0}^r(\hat{x},\hat{y}) = 0 \quad \textrm{ for all } x,y \in
\textup{int } B_r. \] In particular,
\[ p_{t_0}^r(x,x) = \sum_{i=0}^{k(r)} e^{-\lambda_i^r t_0} (\phi_i^r (x))^2 = 0
\quad \textrm{ for all } x \in \textup {int } B_r \] implying that
$\phi_i^r(x) = 0$ for all $i$ and all $x \in$ int $B_r$
contradicting the fact that $ \{\phi_i^r \}_{i=0}^{k(r)}$ forms an
orthonormal basis for $C(B_r, \partial B_r)$.

\bf{6)}  We may assume that $x,y \in$ int $B_r$ since $p_t^r(x,y) =
0$ otherwise. Then, for $x \in $ int $B_r$ and $t=0$, we get
\[ \sum_{y \in \textup{int } B_r} p_0^r(x,y) = \sum_{y \in \textup{int }B_r }
\delta_x(y) = 1. \]  Now, by using the analogue of Green's Theorem
over the interior of $B_r$, whose boundary consists of vertices in
the interior which have a neighbor in the boundary, we will show
that the expression $\sum_{y \in \textup{int } B_r} p_t^r(x,y)$ is
decreasing as a function of $t$. That is,
\begin{eqnarray*}
\pdt \sum_{y \in \textup{int } B_r} p_t^r(x,y) &=& \sum_{y \in
\textup {int } B_r} \pdt p_t^r(x,y) \\
&=& \sum_{y \in \textup{int }B_r} - \Delta_y p_t^r(x,y) \\
&=& \sum_{\substack{y \in \textup {int } B_r \\ z \sim y, \ z \in
\partial B_r}} \Big( p_t^r(x,z) - p_t^r(x,y) \Big) \\
&=& \sum_{\substack{y \in \textup {int } B_r \\ z \sim y, \ z \in
\partial B_r}} - p_t^r(x,y) < 0.
\end{eqnarray*}
Therefore,
\[ \sum_{y \in B_r} p_t^r(x,y) < 1 \quad \textrm{ for all } t>0, \textrm{all } x \in B_r. \]
\qed \newline

\begin{remark}
The proof of Part 6) is identical to the proof of the corresponding
property for the Dirichlet heat kernels on a Riemannian manifold
\cite[Lemma 3.3 Part (i)]{Do2} and shows that a finite graph with
Dirichlet boundary conditions is \emph{not} stochastically complete.
\end{remark}

We now wish to show that the $p_t^r(x,y)$ converge to the heat
kernel $p_t(x,y)$ mentioned at the beginning of this section.  For
this purpose the following lemma will be instrumental.

\begin{lemma}\label{converge}
\[p_t^r(x,y) \leq p_t^{r+1}(x,y) \quad \textrm{ for all } t \geq 0, \textrm{all } x,y
\in B_r.\]
\end{lemma}

\Proof  This is clear for $x$ or $y$ in $\partial
B_r$.  Now, fix $y \in$ int $B_r$ and let
\[ u(x,t) = p_t^{r+1}(x,y) - p_t^r(x,y). \]
Then $\Delta u + \frac{\partial u}{\partial t} = 0$ on int $B_r
\times (0,T)$ which implies that the minimum of $u$ is attained on
the set $(B_r \times \{0\}) \cup (\partial B_r \times [0,T])$.
Since, $u(x,0) = 0$ by Part 4) of Theorem \ref{Dirichlet heat kernels} while, on $\partial B_r \times [0,T]$,
\[ u(x,t) = p_t^{r+1}(x,y) > 0. \]
It follows that
\[ \min_{B_r \times [0,T]} u \geq 0. \]
Therefore, $p_t^{r+1}(x,y) \geq p_t^r(x,y)$ for all $x,y \in B_r$
and for all $t \geq 0$. \qed \newline

By extending $p_t^r(x,y)$ to be 0 outside of $B_r$ and using $0 \leq p_t^r(x,y) \leq 1$ and Lemma \ref{converge} we see that $p_t^r(x,y)$ converges pointwise as $r \to \infty$ for all
$x,y \in V$ and all $t \geq 0$.  Let $p_t(x,y)$ be the limit
\[ p_t^r(x,y) \to p_t(x,y) \textrm{ as } r \to \infty.  \]
We will now show that the convergence is uniform in $t$ for every
compact interval $[0,T]$. To this end, we fix $x$ and $y$ in $V$,
let $f_r(t) = p_t^r(x,y)$ and $f(t) = p_t(x,y)$. Then, from the
definition and properties of each of the heat kernels $p^r$, we get that
each $f_r: [0, \infty) \to \mathbf{R}$ is $C^{\infty}$ and satisfies
\[ \begin{array}{lll}
1) \ f_r \leq f_{r+1} \\
2) \ f_r(t) \to f(t) \textrm{ pointwise for all } t \\
3) \ f_r(t) \leq 1. \\
\end{array}\]
Dini's Theorem implies that $f_r \to f$ uniformly on all compact
subsets $[0,T] \subset [0,\infty)$.

We will now show that $p_t(x,y)$ satisfies the heat equation.
This will follow if we are able to show that $\pdt p_t^r(x,y)$
converges uniformly in $t$ on compact intervals as $r \to \infty$.
But
\begin{eqnarray*}
\pdt p_t^r(x,y) &=& -\Delta_x p_t^r(x,y) \\
&=& \sum_{z \sim x} \big( p_t^r(z,y) - p_t^r(x,y) \big)
\end{eqnarray*}
and since both $p_t^r(z,y)$ and $p_t^r(x,y)$ converge uniformly in
$t$ on $[0,T]$ it follows that $\pdtp_t(x,y)$ exists and is
continuous. In fact, iterating this argument and using the fact that
$\frac{\partial ^ i}{\partial t^i} p_t^r(x,y)$ also satisfy the heat
equation and are continuous for all $i$ we get that $p_t(x,y)$ is
$C^\infty$ in $t$. Then, from the pointwise convergence of
$p_t^r(x,y)$, we get that
\begin{eqnarray*}
\pdt p_t(x,y) &=& \pdt \lim_{r \to \infty} p_t^r(x,y) \\
&=& \lim_{r \to \infty} \pdt p_t^r(x,y) \\
&=& \lim_{r \to \infty} -\Delta_x p_t^r(x,y) \\
&=& -\Delta_x p_t(x,y)
\end{eqnarray*}
implying that $\Delta_x p_t(x,y) + \pdt p_t(x,y) = 0$. The same
argument applied in the $y$ variable then gives $\Delta_y p_t(x,y) +
\pdt p_t(x,y) = 0$.  In summation, using the corresponding
properties of $p_t^r(x,y)$, Lemma \ref{converge}, and what was just
shown, we have proved statements 1), 2), 3), 4), 5), and 6) of the
following theorem. \newline

\begin{theorem}\label{const.heat.kernel}
$p: V \times V \times [0,\infty) \to \mathbf{R}$ has the following
properties:
\begin{enumerate}
\item[\textup{1)}]  $p_t(x,y) > 0$ and $p_t(x,y) = p_t(y,x)$ for all $t > 0$,
all $x,y \in V$.
\item[\textup{2)}]  $p$ is $C^{\infty}$ in $t$.
\item[\textup{3)}] $\Delta p_t(x,y) + \pdt p_t(x,y) = 0$
where $\Delta$ denotes the Laplacian in either $x$ or $y$.
\item[\textup{4)}] $p_0(x,y) = \delta_x(y)$ for all $x,y \in V$.
\item[\textup{5)}]  $p_{s+t}(x,y) = \sum_{z \in V}p_s(x,z)
p_t(z,y)$.
\item[\textup{6)}]  $\sum_{y \in V} p_t(x,y) \leq 1$ for all $t
\geq 0$, all $x,y \in V$.
\item[\textup{7)}] $p$ is independent of the exhaustion used to
define it.
\item[\textup{8)}] $p$ is the smallest non-negative function that satisfies Properties
3) and 4).
\end{enumerate}
\end{theorem}

\Proof
\bf{7)}  Say $D_i$ is another exhaustion of $G$.  That is, each
$D_i$ is a finite and connected subgraph, $D_i \subset D_{i+1}$ for all $i$,
and $G = \bigcup_{i=1}^{\infty} D_i$.  Let $q_t^{D_i}(x,y)$ denote
the Dirichlet heat kernels for this exhaustion and say
that $q_t^{D_i}(x,y) \to q_t(x,y)$.  Then for every $D_i$ there
exists $R$ large enough so that $D_i \subset B_R$.  By the maximum
principle, since $q_t^{D_i}(x,y)$ vanishes on $\partial D_i$, we
obtain $q_t^{D_i}(x,y) \leq p_t^R(x,y)$.  Because $p_t^R(x,y) \leq
p_t^{R+1}(x,y) \leq \dots$ and $p_t^R(x,y) \to p_t(x,y)$ this
implies that $q_t^{D_i}(x,y) \leq p_t(x,y)$.  Letting $i \to \infty$
gives
\[ q_t(x,y) \leq p_t(x,y). \]  Interchanging the roles of
$q^{D_i}$ and $p^r$ in the preceding argument gives $q_t(x,y) \geq
p_t(x,y)$ and therefore, $p_t(x,y) = q_t(x,y)$. \newline

\bf{8)}  Say $q_t(x,y)$ is another non-negative function that
satisfies Properties 3) and 4).  In particular, both $q$ and $p^r$
satisfy the heat equation on int $B_r \times (0,T)$.  Since $p^r$
vanishes on $\partial B_r$ while $q$ is non-negative there we get,
by applying the maximum principle to the difference of $q$ and
$p^r$, that $q_t(x,y) \geq p_t^r(x,y)$ on $B_r \times [0,T]$ and
hence for all $x, y \in V$ and all $t>0$ . Letting $r \to \infty$ we
get that $q_t(x,y) \geq p_t(x,y)$. \qed
\newline


\section{The Spectral Theorem Construction}
As mentioned previously, an alternative way of obtaining the heat
kernels $p_t^r(x,y)$ on $B_r$ with Dirichlet boundary conditions is
through the convergent power series
\[ e^{-t \Delta_r} = I - t\Delta_r + \frac{t^2}{2}\Delta_r^2 -
\frac{t^3}{6}\Delta_r^3 + \ldots \] by letting $p_t^r(x,y) = (e^{-t
\Delta_r} \delta_x)(y).$  However, on the entire graph, since
the Laplacian is not bounded on $\ell^2(V)$, one cannot use the power
series approach. One can still construct $e^{-t \tilde{\Delta}}$ for, $\tilde{\Delta}$, the
unique self-adjoint extension of $\Delta$ to $\ell^2(V)$ by using
the functional calculus developed through the spectral theorem
\cite[Chapter VIII]{Reed-Simon}.  The purpose of this section is to
show that the construction given in the previous section via exhaustion and this
approach result in the same kernel \cite[Proposition 4.5]{Do2}. Let
$P_t v(x) = \sum_{y \in V} p_t(x,y) v(y)$ for any bounded function
$v$, with $P_t^r v$ indicating a similar sum for the heat kernel
$p_t^r(x,y)$ on $B_r$. We then have the following theorem which
states that $P_t$ and $e^{-t \tilde{\Delta}}$  agree on a dense subset of
$\ell^2(V)$ and, as such, have the same kernel:

\begin{theorem}
\[ P_t v = e^{-t \tilde{\Delta}} v \quad \textrm{for all } v \in
C_0(V).\]
\end{theorem}
\Proof We begin by showing that if $v \in C_0(V)$ then $P_t v \in
\ell^2(V)$ and $\Delta P_t v \in \ell^2(V)$. Since $v$ is finitely
supported, there exists a ball of large radius $R$ which contains
its support. Therefore,
\begin{eqnarray*}
\| P_t^R v \|^2_{\ell^2(V)} &=& \sum_{x \in V} (P_t^R v (x))^2 \\
&=& \sum_{x \in V} \left( \sum_{y \in B_R} p_t^R(x,y) v(y) \right)^2
\\
&\leq& \sum_{y \in B_R} \left( \sum_{x \in B_R} p_t^R(x,y)^2 \right)
v(y)^2 \\
&\leq& \sum_{y \in B_R} v(y)^2 = \| v \|^2_{\ell^2(V)}.
\end{eqnarray*}
By letting $R \to \infty$ and using the dominated convergence
theorem it follows that $P_t v \in \ell^2(V).$  In fact, this actually proves that $P_t$ is a bounded operator on $\ell^2(V)$ with $\| P_t \| \leq 1.$

We will now show that $\Delta P_t v = P_t \Delta v$ and since if $v$ is finitely supported then so is $\Delta v$ it will follows that $\Delta P_t v \in \ell^2(V)$. To show that $\Delta P_t
v = P_t \Delta v$ we calculate:
\begin{eqnarray*}
\Delta(P_t v)(x) &=& \sum_{y \sim x} \big( (P_t v)(x) - (P_t v)(y)
\big) \\
&=& \sum_{y \sim x} \left( \sum_{z \in V}\big(p_t(x,z) - p_t(y,z)\big) v(z)
\right).
\end{eqnarray*}  Meanwhile, by using the analogue of Green's Theorem and the fact
that the heat kernel satisfies the heat equation in both variables,
we get
\begin{eqnarray*}
P_t(\Delta v)(x) &=& \sum_{z \in V} p_t(x,z) \Delta v(z) \\
&=& \sum_{z \in V} \Delta_z p_t(x,z) v(z) \\
&=& \sum_{z \in V} \Delta_x p_t(x,z) v(z) \\
&=& \sum_{z \in V} \left( \sum_{y \sim x} \big(p_t(x,z) - p_t(y,z) \big)
\right) v(z).
\end{eqnarray*}

We now give the proof of the theorem.  Let
\[ u(x,t) = \left( P_t - e^{-t \tilde{\Delta}} \right) v(x). \]
Then $u(x,0) = 0$ and since $\Delta P_t v$ is in $\ell^2(V)$ we can
apply Green's Theorem again to obtain
\begin{eqnarray*}
\pdt \sum_{x \in V} u^2(x,t) &=& 2 \sum_{x \in V}u(x,t)
\pdt u(x,t) \\
&=& -2\sum_{x \in V} u(x,t) \Delta u(x,t) \\
&=& - 2\sum_{[x,y] \in \tilde{E}} \big( u(y,t) - u(x,t) \big)^2 \leq 0
\end{eqnarray*}
from which it now follows that $P_t v(x) = e^{-t \tilde{\Delta}} v(x)$ for
all finitely supported $v$.  Since both $P_t$ and $e^{-t \tilde{\Delta}}$ are bounded it follows that they are equal on $\ell^2(V)$.  \qed


\chapter{Stochastic Incompleteness}

\section{Stochastic Incompleteness}

We now define the notion of stochastic incompleteness and
recall the proof of the equivalence of several properties and this definition.
The material here is adapted from \cite[$\textrm{p.}$ 170-172]{Grig}. We recall
the definition of $P_t$:
\[ P_t u_0(x) = \sum_{y \in V} p_t(x,y) u_0(y) \] for any bounded
function $u_0$ on $G$.  This summation converges from Part 6) of Theorem \ref{const.heat.kernel} and  from Part 5), $P_t$ satisfies the semigroup property:
\[ P_s (P_t u_0) = P_{s+t} u_0. \]
Apply $P_t$ to the function \bf{1} which is exactly 1 on each
vertex of $G$:
\[ P_t \bf{1}(x) = \sum_{y \in V} p_t(x,y) \]
and note that this sum is less than or equal to 1 from Part 6) of Theorem
\ref{const.heat.kernel}.

\begin{definition}
A graph $G$ is called \emph{stochastically incomplete} if for some
vertex $x_0$ of $G$ and some $t_0>0$
\[ P_{t_0} \bf{1}(x_0) = \sum_{y \in V} p_{t_0}(x_0,y) < 1. \]
\end{definition}

\begin{remark}
Although this really is a property of the heat kernel or of the diffusion process which is
modeled by the heat kernel, it is customary to say, as above, that
it is a property of the underlying space.
\end{remark}

\begin{theorem}\label{stochastic.incomp}
The following statements are equivalent:
\begin{enumerate}
\item[\up{1)}] For some $t_0>0,$ some $x_0 \in V$, $P_{t_0} \up{\bf{1}}(x_0) < 1. $
\item[\up{1')}] For all $t>0$, all $x \in V$, $P_t \up{\bf{1}}(x)
< 1. $
\item[\up{2)}] There exists a positive (equivalently, non-zero) bounded function $v$ on
$G$ such that $\Delta v = \lambda v$ for any $\lambda < 0$.
\item[\up{2')}] There exists a positive, (equivalently, non-zero) bounded function $v$ on
$G$ such that $\Delta v \leq \lambda v$ for any $\lambda < 0$.
\item[\up{3)}]  There exists a nonzero, bounded solution to
\[ \left\{  \begin{array}{ll} \Delta u(x,t) + \pdtu(x,t) = 0 & \textrm{ for all } x \in V,
\textrm{all } t>0\\
 u(x,0) = 0 & \textrm{ for all }  x \in V.
 \end{array} \right.  \] \\
\end{enumerate}
\end{theorem}

\begin{definition}
Any function $v$ on $G$ such that $\Delta v = \lambda v$ is called
$\lambda$\emph{-harmonic} whereas if $\Delta v \leq \lambda v$, $v$
is called $\lambda$\emph{-subharmonic}.
\end{definition}

Therefore, stochastic incompleteness is equivalent to the existence of a positive,
bounded $\lambda$-harmonic (or $\lambda$-subharmonic) function for negative $\lambda$ and to the non-uniqueness of bounded solutions for the heat equation on $G$.

\Proof
\newline \noindent \bf{1') $\Rightarrow$ 1)}  Obvious.
\newline \noindent \bf{1) $\Rightarrow$ 1')}  If there exists $x_0 \in
V$ and a $t_0 > 0$ such that $P_{t_0}\bf{1}(x_0) = 1$ then by the stong
maximum principle for the heat equation, Lemma \ref{maxheatlem}, applied to the function $P_t
\bf{1}$ we get that
\[ P_{t_0} \bf{1}(x) = 1 \textrm{ for all } x.\]

Now, if $s < t_0$ then it follows from the semigroup property that
\[ P_{t_0} \bf{1} = P_{t_0 - s} (P_s \bf{1}) \leq P_{t_0 - s} \bf{1} \leq 1. \]
For any $t_0 > 0$ such that $P_{t_0} \bf{1} = 1$ it
follows that the inequalities become equalities and, in particular,
$P_s \bf{1} = 1$ for all $s < t_0$. If $s > t_0$ then there exists a
$k$ such that $s < kt_0$ and by the semigroup property we get
that
\[ P_{kt_0}\bf{1} = (\underbrace{P_{t_0} \ldots P_{t_0}}_k)\bf{1} =1 \]
provided that $P_{t_0}\bf{1} = 1$ giving $P_s \bf{1} = 1$ from the
same argument as above. \newline

\noindent \bf{1') $\Rightarrow$ 2)}  For any $\lambda < 0$, let
$w(x) = \int_0^{\infty} e ^ {\lambda t} u(x,t) dt \ $ where $u(x,t)
= P_t \bf{1}(x) < 1$ by assumption.  Then
\begin{eqnarray*}
0 < w &<& \int_0^{\infty} e ^ {\lambda t} dt\\
&=& \frac{1}{\lambda}\left(e^{\lambda t} \Big|_0^{\infty}\right) \\
&=& \frac{1}{\lambda}(0 - 1) = -\frac {1}{\lambda}.
\end{eqnarray*}
Integration by parts gives
\begin{eqnarray*}
\Delta w = \int_0^{\infty} e ^ {\lambda t} \Delta u(x,t) dt &=&
-\int_0^{\infty} e ^ {\lambda t} \pdtu(x,t) dt\\
&=& -e^{\lambda t}u(x,t)\Big|_0^{\infty} + \int_0^{\infty}\lambda e ^ {\lambda t}u(x,t) dt\\
&=& 1 + \lambda w.
\end{eqnarray*}
If $v=1+\lambda w$, then $v$ satisfies
\[ \Delta v = \lambda \Delta w = \lambda(1 + \lambda w) = \lambda
v \] which shows that $v$ is $\lambda$-harmonic.  Since $0 < w <
-\frac{1}{\lambda}$ we have $0 < v < 1$ so that $v$ is positive and
bounded. \newline

\noindent \bf{2) $\Rightarrow$ 2')}  Clear.
\newline \noindent \bf{2') $\Rightarrow$ 2)}  Exhaust the graph $G$ by
finite, connected subgraphs $D_i$.  That is, $D_i \subset D_{i+1}$
and $G = \bigcup_{i=0}^{\infty} D_i$ where each $D_i$ is finite and
connected.  Let $\Delta_i$ denote the reduced Laplacian acting on
the space $C(D_i,\partial D_i)$ of functions on $D_i$ which vanish
on the boundary $\partial D_i$. Then, for $\lambda < 0$, one can
solve
\begin{equation}\label{lambda.subhar}
\left\{  \begin{array}{ll} \Delta_i v_i = \lambda v_i & \textrm{on int } D_i\\
 {v_i}_{| \partial D_i} = 1.
 \end{array} \right.
\end{equation}
Indeed, letting $\bf{1}_{D_i}$ denote the function that is 1 on
every vertex of $D_i$ and 0 elsewhere, if $v_i$ is a solution to the
above then $w_i = v_i - \bf{1}_{D_i}$ would vanish on the boundary
of $D_i$ and on the interior would satisfy
\[ \Delta_i w_i(x) = \Delta_i v_i(x) = \lambda v_i(x) = \lambda(w_i(x) + 1). \]
That is,
\[ (\Delta_i - \lambda I) w_i = \lambda_{\textrm{int } D_i} \]
where $\lambda_{\textrm{int } D_i}$  denotes the function that is equal to
$\lambda$ on every vertex in the interior of $D_i$ and is $0$ on
$\partial D_i$.  Since $\lambda < 0$, $\Delta_i - \lambda I$ is
invertible on $C(D_i,\partial D_i)$ and so
\[ w_i = (\Delta_i - \lambda I)^{-1} (\lambda_{\textrm{int } D_i}) \]
yielding
\[ v_i = (\Delta_i - \lambda I)^{-1} (\lambda_{\textrm{int } D_i}) + \bf{1}_{D_i} \]
as a solution for (\ref{lambda.subhar}).

We now claim that \begin{equation} \label{v_i leq 1}
0 < v_i \leq 1 \textrm{ on } D_i.
\end{equation}
This follows from the
fact that, if there exists an $x_0$ in the interior of $D_i$ such
that $v_i(x_0) \leq 0$, then we may assume that $x_0$ is a minimum
for $v_i$ and
\[ \Delta v_i(x_0) = \sum_{x \sim x_0} \big( v_i(x_0) - v_i(x) \big)
\leq 0 \] while $\Delta v_i(x_0) = \lambda v_i(x_0) \geq 0$ so that
$\Delta v_i(x_0)=0$.  This implies that $v_i(x) = v_i(x_0)$ for all
neighbors of $x_0$ and by repeating the argument we would get that
$v_i$ is a non-positive constant on $D_i$ contradicting that $v_i =1 $ on the boundary of
$D_i$.

Therefore, $v_i > 0$ and so $\Delta v_i < 0$ on the interior which
implies that
\[ \max_{D_i} v_i = \max_{\partial D_i} v_i = 1. \]  Indeed, at an
interior maximum $\Delta v_i(x_0) \geq 0$.  This completes the proof of
(\ref{v_i leq 1}).

Furthermore, if we extend each $v_i$ to be exactly 1 outside of $D_i$, it is true that
\[v_i \geq v_{i+1}.\]
This is clear on $\partial D_i$ since $v_i =1$ there while $v_{i+1}
\leq 1.$ On the interior of $D_i$ we have that $\Delta (v_i -
v_{i+1}) = \lambda (v_i - v_{i+1})$ from which it follows that $v_i
- v_{i+1} > 0$ by the same argument that gives $v_i > 0$ above.

The $v_i$ therefore, form a non-increasing, bounded sequence so that
\[v_i \to v\]
where $0 \leq v \leq 1$ and $\Delta v = \lambda v$. What remains to
be shown is that $v$ is positive (or non-zero) and here we use the assumption that there
exists on $G$ a positive (or non-zero), bounded $\lambda$-subharmonic function
$w$. That is, $w$ is positive (or non-zero), bounded and satisfies $\Delta w \leq
\lambda w$ on $G$.  Assuming that $w \leq 1$, we show that
\[  v_i \geq w \textrm{ on } D_i \textrm{ for all } i. \]
This can be seen as follows: on the interior of $D_i$,
\begin{equation} \label{v_i and w}
\Delta (v_i - w) \geq \lambda (v_i - w).
\end{equation}
Therefore, if there exists an $x_0$ in the interior of $D_i$ such that $(v_i - w)(x_0) < 0$ and $x_0$ is
a minimum for $v_i - w$ then by computation $\Delta (v_i - w)(x_0)
\leq 0$ while $\Delta(v_i - w)(x_0) \geq \lambda(v_i - w)(x_0) > 0$
from (\ref{v_i and w}). The contradiction implies that $v_i \geq
w$ on $D_i$ and by passing to the limit we get that $v \geq w$ so
that $v$ is either positive or non-zero depending on $w$.
\newline

\noindent \bf{2)$\Rightarrow$ 3)}  Let $w(x,t) = e^{-\lambda t}
v(x)$ where $v$ is a positive, bounded $\lambda$-harmonic function
for $\lambda < 0$. Then $w$ is positive and bounded on $V \times [0,
T]$, and satisfies $ \Delta w + \pdtw = 0 $ with $w(x,0) = v(x).$
The function $P_t v$ also satisfies both equations and,
moreover,
\begin{equation} \label{inequality1}
\sup_{x \in V} P_t v(x) \leq \sup_{x \in V} v(x)
\end{equation}
while
\begin{equation} \label{inequality2}
w(x,t) > v(x) \textrm{ for all } t > 0, \textrm{ all } x \in V.
\end{equation}
Therefore, we have two different bounded solutions to
\[ \left\{  \begin{array}{ll} \Delta u(x,t) + \pdtu(x,t) = 0
& \textrm{for all } (x,t) \in V \times (0,T)\\
 u(x,0) = v(x) & \textrm{all } x \in V.
 \end{array} \right. \]
Taking the difference of the two solutions gives a nonzero, bounded
solution to
\begin{equation}\label{heat0}
\left\{  \begin{array}{ll} \Delta u(x,t) + \pdtu(x,t) = 0
& \textrm{for all } (x,t) \in V \times (0,T)\\
 u(x,0) = 0 & \textrm{all } x \in V.
 \end{array} \right.
\end{equation}

Therefore, we have shown the existence of a nonzero, bounded
solution for the heat equation with initial condition 0 for a finite
time interval. In the argument below, we show that this is enough to
imply condition 1), that is, $P_t \bf{1} < 1$. Then, given $P_t
\bf{1} < 1$, $\bf{1} - P_t \bf{1}$ will give a nonzero,
bounded solution to (\ref{heat0}) for an infinite time interval,
completing the proof.

We also note that the assumption that $v$ is positive was not
essential for the argument as can be seen by putting norms about
$P_t v$, $v$, and $w$ in (\ref{inequality1}) and (\ref{inequality2}).
Therefore, the existence of any bounded, \emph{non-zero}
$\lambda$-harmonic function will imply stochastic incompleteness and
then the argument giving the implication $1') \Rightarrow 2)$ shows
that there then exists a \emph{positive}, bounded $\lambda$-harmonic
function on $G$. \newline

\noindent \bf{3) $\Rightarrow$ 1)}  Suppose that $u(x,t)$ is nonzero,
bounded and satisfies
\[ \left\{  \begin{array}{ll} \Delta u(x,t) + \pdtu(x,t) = 0 &
\textrm{for all }(x,t) \in V \times (0,T)\\
 u(x,0) = 0 & \textrm{all } x \in V.
 \end{array} \right.\]
Then, by rescaling, we may assume that $|u(x,t)| < 1$ for all $x$ and $t$, and
that there exists an $x_0$ and $t_0 > 0$ such that $u(x_0,t_0) > 0$.
Then, $w(x,t) = 1 - u(x,t)$ is bounded,
positive and satisfies
\begin{equation} \label{heat1}
\left\{  \begin{array}{ll} \Delta w(x,t) + \pdtw(x,t) = 0 & \textrm{for all }(x,t) \in V \times (0,T)\\
 w(x,0) = 1 & \textrm {all } x \in V.
 \end{array} \right.
\end{equation}
Furthermore, $w(x_0,t_0) < 1$.  Now, letting $P_t^r \bf{1}(x) =
\sum_{y \in B_r}p_t^r(x,y)$, and applying the maximum principle for
the heat equation to $P_t^r \bf{1}(x) - w(x,t)$ on $B_r \times
[0,T]$, we get that $P_t^r \bf{1}(x) < w(x,t)$ for all $r$.
Therefore, letting $r \to \infty$,
\[ P_{t_0} \bf{1}(x_0) \leq w(x_0,t_0) < 1 .\] \qed


\section{Model Trees}

We now turn our focus to a family of particular trees and study
under what conditions they are stochastically complete. A tree will
be called \emph{model} if it contains a vertex $x_0$, henceforth
called the \emph{root} for the model, such that the valence $m(x)$
is constant on spheres $S_r(x_0) = S_r$ of radius $r$ about $x_0$. That is, if
\[ S_r = \{ x \ | \ d(x,x_0) = r \} \] then
\[ m(x) = m(r) \textrm{ for all } x \in S_r. \]
For $r>0$, we let $n(r) = m(r) - 1 $ denote the \emph{branching} of
$T$, that is, the number of edges connecting a vertex in $S_r$ with
vertices in $S_{r+1}$, and let $n(0) = m(x_0)$. We then denote such
trees $T_n$ to indicate that their structure is completely encoded
in the branching function $n(r)$.

The main purpose of this section is to prove the following theorem which
tells us precisely when such trees are stochastically complete.

\begin{theorem}\label{models}
$T_n$ is stochastically complete if and only if
\[ \sum_{r=0}^{\infty} \frac{1}{n(r)} = \infty. \]
\end{theorem}

The idea of the proof is to study positive $\lambda$-harmonic
functions on $T_n$ for $\lambda < 0$.  By averaging over spheres we
can reduce to the case of $\lambda$-harmonic functions depending
only on the distance $r$ from the root $x_0$.  It then turns out that
such a function will be bounded if and only if the series above
converges.  Since the existence of a positive, bounded, $\lambda$-harmonic
function is equivalent to stochastic completeness, Theorem \ref{models} will
then follow.

Let, therefore, $v(r)$ denote a function on the vertices of $T_n$
depending only on $r = r(x) = d(x,x_0)$, the distance from a vertex
to the root.  When the Laplacian is applied to such a function we get
\begin{eqnarray*}
\Delta v(r) &=& \big( n(r)+1 \big) v(r) - n(r) v(r+1) - v(r-1) \\
&=& n(r) \Big(v(r) - v(r+1) \Big) + \Big(v(r) - v(r-1) \Big).
\end{eqnarray*}
We will study the existence and boundedness of such functions on
$T_n$ when, in addition, they are positive and $\lambda$-harmonic
for a negative $\lambda$, that is, they satisfy $\Delta v(r) = \lambda v(r)$ for
$\lambda < 0$. We start by showing that there is no loss in
generality in restricting our study of $\lambda$-harmonic functions to
only functions of this type.

\begin{lemma}\label{averagelemma}
If there exists a positive, bounded $\lambda$-harmonic function on
$T_n$ then there exists one depending only on the distance from
$x_0$.
\end{lemma}
\Proof  Let $u(x)$ denote a positive, bounded $\lambda$-harmonic
function on $T_n$.  If $S_r$ denotes the sphere of radius $r$ about
the root $x_0$ and Vol($S_r$) denotes its volume, that is,
\[ \textrm{Vol}(S_r) = \# \{x \ | \ x \in S_r \} \]
then we define a function $v(r)$
depending only on the radius $r$ by averaging $u$ over $S_r$:
\[ v(r) = \frac{1}{\textrm{Vol}(S_r)} \sum_{x \in S_r} u(x). \]
Clearly, such a function will be positive and bounded since $u(x)$ is.
We now show that $v(r)$ is also $\lambda$-harmonic.  At $x_0$, since Vol$(S_1) = n(0)$, we
have:
\begin{eqnarray*}
\Delta v(0) &=& n(0) \Big( v(0) - v(1) \Big) \\
&=& n(0) \left( u(x_0) - \frac{1}{n(0)} \sum_{x \in S_1} u(x) \right) \\
&=& n(0)u(x_0) - \sum_{x \sim x_0} u(x) \\
&=& \Delta u(x_0) = \lambda u(x_0) = \lambda v(0).
\end{eqnarray*}

Now, if $x \in S_r$ for $r>0$, then
\begin{eqnarray*}
\Delta u(x) &=& \big( n(r)+1 \big)u(x) - \sum_{\substack{ z \sim x \\
z \in S_{r+1} } } u(z) - u(y) \\
 &=& \lambda u(x)
\end{eqnarray*} where $y$ is the unique neighbor of $x$ that
is in $S_{r-1}$.  Therefore,
\begin{equation}\label{average}
\big( n(r) + 1 - \lambda \big) u(x) = \sum_{\substack{ z \sim x \\ z
\in S_{r+1} } } u(z) + u(y).
\end{equation}
We will average this equation over $S_r$ and use the following
version of the formula $n(r) \cdot \textrm{Vol}(S_r) = \textrm{Vol}
(S_{r+1})$ which relates the volume of spheres of different radius in
$T_n$:
\begin{equation}\label{volumes}
\frac{1}{\textrm{Vol}(S_r)} = \frac{n(r)}{\textrm{Vol}(S_{r+1})}.
\end{equation}

From the definition of $v$, we get
\begin{eqnarray*}
\big( n(r) + 1 - \lambda \big) v(r) &=& \frac{(n(r) + 1 -
\lambda)}{\textrm{Vol}(S_r)} \sum_{x\in S_r} u(x) \\
&=& \frac{1}{\textrm{Vol}(S_r)} \sum_{x \in S_r} \left(
\sum_{\substack{ z \sim x \\ z \in S_{r+1} } } u(z) +
\sum_{\substack{ y \sim x \\ y \in S_{r-1}}} u(y) \right) \quad \ \ \ \ \ (\ref{average}) \\
&=& \frac{1}{\textrm{Vol}(S_r)} \sum_{z \in S_{r+1}} u(z) +
\frac{n(r-1)}{\textrm{Vol}(S_r)} \sum_{y \in
S_{r-1}} u(y) \\
&=& \frac{n(r)}{\textrm{Vol}(S_{r+1})} \sum_{z \in S_{r+1}} u(z) +
\frac{1}{\textrm{Vol}(S_{r-1})} \sum_{y \in S_{r-1}} u(y) \quad (\ref{volumes})\\
&=& n(r)v(r+1) + v(r-1)
\end{eqnarray*}
or precisely that $\Delta v(r) = \lambda v(r)$. \qed \newline

Therefore, on a model tree, the existence of any positive, bounded $\lambda$-harmonic
function is equivalent to the existence of such a function
depending only on the distance to the root.  The values of such a function are determined by the value of the function at the root and are given by:
\begin{eqnarray}
v(1) &=& \left( 1 - \frac{\lambda}{n(0)} \right) v(0) \label{v(0)} \\
v(r+1) &=& \frac{1}{n(r)} \Bigg( \big(n(r)+1-\lambda \big)v(r) -
v(r-1) \Bigg). \label{v(1)}
\end{eqnarray}
We will now study under what conditions such a function
will remain bounded. We start by showing that such a function must
increase with the radius.

\begin{lemma}\label{increasing}
If $v > 0$ satisfies $\Delta v = \lambda v$ for $\lambda < 0$ then
\[v(r) < v(r+1) \textrm{ for all } r \geq 0.\]
\end{lemma}
\Proof  The proof is by induction.  We have that $v(0) < v(1)$ from
(\ref{v(0)}).

Now, assuming that $v(r-1) < v(r)$,
\[ \Delta v(r) = n(r) \Big(v(r) - v(r+1) \Big) + \Big(v(r) - v(r-1)
\Big) = \lambda v(r) \] gives
\[ n(r) \Big( v(r) - v(r+1) \Big) = \lambda v(r) - \Big( v(r) -
v(r-1) \Big) < 0 \] implying
\[ v(r) < v(r+1). \] \qed

\begin{lemma}\label{product}
If $\Delta v = \lambda v$ with $v>0$ and $\lambda < 0$ then
\[ \prod_{i=0}^r \left( 1 - \frac{\lambda}{n(i)} \right) v(0) < v(r+1) <
\prod_{i=0}^{\infty} \left( 1+ \frac{1 - \lambda}{n(i)} \right)
v(0). \]
\end{lemma}

Consequently, for $\lambda$ negative, a positive,
$\lambda$-harmonic function on $T_n$ depending only on the distance
from the root remains bounded if and only if $\
\prod_{i=0}^{\infty}\left( 1 + \frac{1}{n(i)} \right) < \infty$.
That is, the following conditions are equivalent:
\begin{enumerate}
\item[1)] $v(r)$ is bounded
\item[2)] $\prod_{i=0}^{\infty}\left( 1 + \frac{1}{n(i)} \right) < \infty$
\item[3)] $ \sum_{i=0}^{\infty}\frac{1}{n(i)} < \infty.$
\end{enumerate}

\Proof  For the upper bound, we rewrite the relation $\Delta v(r) =
\lambda v(r)$ as:
\[ \Big( n(r) + 1 - \lambda \Big) v(r) - n(r) v(r+1) = v(r-1) > 0. \]
Therefore,
\[ \Big( n(r) + 1 - \lambda \Big) v(r) > n(r) v(r+1) \]
or
\begin{eqnarray*}
v(r+1) &<& \frac{(n(r) + 1 - \lambda)}{n(r)} v(r) \\
&=& \left( 1 + \frac{1-\lambda}{n(r)} \right) v(r).
\end{eqnarray*}
Now, iterate this relation down to $v(0)$:
\begin{eqnarray*}
v(r+1) &<& \left( 1 + \frac{1-\lambda}{n(r)} \right) v(r) \\
&<& \left( 1 + \frac{1 - \lambda}{n(r)} \right) \left( 1 +
\frac{1-\lambda}{n(r-1)} \right) v(r-1) \\
&<& \prod_{i=0}^{r} \left( 1+ \frac{1 - \lambda}{n(i)} \right) v(0) \\
&<& \prod_{i=0}^{\infty} \left( 1+ \frac{1 - \lambda}{n(i)} \right)
v(0).
\end{eqnarray*}

For the lower bound, we use Lemma \ref{increasing}, which implies
that \\ $v(r) - v(r-1) > 0$ as follows:
\begin{eqnarray*}
\Delta v(r) &=& n(r) \Big( v(r) - v(r+1) \Big) +
\Big( v(r) - v(r-1) \Big) \\
&>& n(r) \Big( v(r) - v(r+1) \Big).
\end{eqnarray*}
Since $\Delta v(r) = \lambda v(r)$, this gives
\[ n(r) \Big( v(r) - v(r+1) \Big) < \lambda v(r)\]
or
\[ \left( 1 - \frac{\lambda}{n(r)} \right) v(r) < v(r+1). \]
Iterating as before gives
\[ \prod_{i=0}^r \left( 1 - \frac{\lambda}{n(i)} \right) v(0) < v(r+1) \]
completing the proof of the lemma. \qed \newline

\noindent \textbf{Proof of Theorem \ref{models}:} By Theorem
\ref{stochastic.incomp}, stochastic incompleteness is equivalent to
the existence of a positive, bounded $\lambda$-harmonic function for
$\lambda < 0$. We can define such a function on $T_n$ depending
only on the distance from $x_0$ by (\ref{v(0)}) and
(\ref{v(1)}). If $\sum_{r=0}^{\infty}\frac{1}{n(r)} < \infty$, this
function will remain bounded by Lemma \ref{product}.

Now, if $\sum_{r=0}^{\infty} \frac{1}{n(r)} = \infty$ then every
positive, $\lambda$-harmonic function depending only on the radius
from the root will be unbounded.  Therefore, every positive,
$\lambda$-harmonic function on $T_n$ will be unbounded by Lemma
\ref{averagelemma} so that $T_n$ is stochastically complete.  \qed
\\

\begin{remark}
We would like to point out the relationship between Theorem \ref{models} and the case of \emph{spherically
symmetric} or \emph{model} manifolds on which we base our definition
of model trees. $M_\sigma$, a Riemannian manifold of dimension $d$
with pole $o$, is called \emph{model} if
\begin{enumerate}
\item[i)]  topologically, $M_\sigma \setminus \{o\}$ is the product of an
open interval $I$ and the sphere $S^{d-1}$.  Therefore, each point
$x \in M_\sigma \setminus \{o\}$ can be identified with a pair
$(r,\theta)$ where $r \in I$ and $\theta \in S^{d-1}$.
\item[ii)] the metric on $M_\sigma$ is given by
\begin{equation}\label{modelmetric}
ds^2 = dr^2 + \sigma^2(r) d\theta^2
\end{equation}
where $d\theta^2$ denotes the standard Euclidean metric on
$S^{d-1}$. Here, $\sigma$ is a smooth, positive function on $I$
sometimes called the \emph{twisting} or \emph{warping function}
\cite[$\textrm{p.}$ 145-148]{Grig}.
\end{enumerate}

It follows from (\ref{modelmetric}) that the area of a sphere of
radius $r$ in $M_\sigma$ is given by
\begin{equation}\label{modelspheres}
A(S_r) = \omega_d \sigma^{d-1}(r)
\end{equation}
where $\omega_d$ is the area of the unit sphere in $\mathbf{R}^d$.
Now, it is shown in \cite[Corollary 6.8]{Grig} that a geodesically
complete, noncompact, model manifold is stochastically complete if
and only if
\[ \int^{\infty} \frac{\textrm{Vol}(B_r)}{A(S_r)} dr = \infty \]
where Vol$(B_r)$ denotes the Riemannian volume of the geodesic ball in
$M_\sigma$.  For example, if, for large $r$, $A(S_r) \leq e^{r^2}$ then
$M_\sigma$ is stochastically complete, whereas, if $A(S_r) =
e^{r^{2+\epsilon}}$ for any positive $\epsilon$, then $M_\sigma$
will be incomplete.

Observe from (\ref{modelspheres}) that, on $M_\sigma$,
\begin{eqnarray*}
dA(S_r) &=& \omega_d (d-1) \sigma^{d-2}(r) \sigma^{\prime}(r) \ dr \\
&=& (d-1) \omega_d (\ln \sigma(r))^\prime A(S_r) \ dr
\end{eqnarray*}
so that
\begin{equation}\label{model.manifold.spheres}
\frac{dA(S_r)}{A(S_r)} = c(\ln \sigma(r))^\prime \ dr
\end{equation}
for $c = (d-1) \ \omega_d$.

Meanwhile, for model trees, temporarily using the notation $A(S_r) =
\textrm{Vol}(S_r)$ and $dA(S_r) = A(S_{r+1}) - A(S_r)$, it follows that
\[ dA(S_r) = \big(n(r) - 1\big)A(S_r) \] so that
\begin{equation}\label{model.tree.spheres}
\frac{dA(S_r)}{A(S_r)} = \big( n(r) - 1 \big).
\end{equation}
Therefore, comparing (\ref{model.manifold.spheres}) and (\ref{model.tree.spheres}) we see
that $n(r)$ and $(\ln \sigma(r))^\prime$ play a similar role and the correspondence between the borderlines for stochastic completeness is exact.

\end{remark}


\section{Comparison Theorems}
Throughout this section, we assume that $T_n$ denotes a model tree
with root vertex $x_0$ while $T$ denotes a general tree.  From Theorem \ref{models} in the
last section, we know
that $T_n$ will be stochastically complete if and only if
$\sum_{r=0}^\infty \frac{1}{n(r)} = \infty$ and we wish to obtain a
similar criterion for $T$.  In order to state our results for
$T$ we make the following definitions.

\begin{definition} For a vertex $x_0 \in T$, let
\begin{eqnarray*}
\underline{m}(r) = \underline{m}_{x_0}(r) &=& \min_{x \in S_r(x_0)}
m(x) \\
M(r) = M_{x_0}(r) &=& \max_{x \in S_r(x_0)} m(x).
\end{eqnarray*}
\end{definition}

The following result is an immediate consequence of our criterion for model
trees and the characterization of stochastic incompleteness in terms
of $\lambda$-subharmonic functions.

\begin{theorem}\label{comp1} Assume that $T_n \subseteq T$ and that $\underline{m}_{x_0}(r) = \underline{m}(r) = \min_{x \in S_r(x_0)} m(x)$ satisfies
\[ n(r) \leq \ \underline{m}(r) - 1 \textrm{ for all } r >0. \]
Then, if $T_n$ is stochastically incomplete, so is $T$.
\end{theorem}
\Proof Since $T_n$ is stochastically incomplete, there exists a
bounded, positive function $v(r)$ on $T_n$ such that $v(r) < v(r+1)$ and $\Delta v(r) = \lambda v(r)$ for $\lambda <
0$. Let $r(x) = d(x,x_0)$ be the distance between $x_0$ and $x \in
T$ and define a function $u$ on $T$ by
\[u(x) = v(r(x)).\]
Clearly, $u(x)$ will be bounded and positive since $v$ is. Now, it
follows from the inequalities $v(0) - v(1) < 0$ and $n(0) \leq
m(x_0)$ that $u(x)$ is $\lambda$-subharmonic at $x=x_0$:
\begin{eqnarray*}
\Delta u(x_0) &=& m(x_0) u(x_0) - \sum_{x \sim x_0} u(x) \\
&=& m(x_0) \Big( v(0) - v(1) \Big) \\
&\leq& n(0) \Big( v(0) - v(1) \Big) \\
&=& \lambda v(0) = \lambda u(x_0).
\end{eqnarray*}

Now, suppose that $r(x) = r > 0$ and $y$ denotes the unique
neighbor of $x$ in $S_{r-1}$.  Then, since $n(r) \leq m(x) - 1$,
\begin{eqnarray*}
\Delta u(x) &=& m(x)u(x) -
\sum_{\substack{z \sim x \\ z \in S_{r+1}}}u(z) - u(y) \\
&=& (m(x)-1) \Big( v(r) - v(r+1) \Big) + \Big(v(r) - v(r-1) \Big) \\
&\leq& n(r) \Big( v(r) - v(r+1) \Big) + \Big( v(r) - v(r-1) \Big) \\
&=& \lambda v(r) = \lambda u(x).
\end{eqnarray*}
Thus, $u$ is a positive, bounded $\lambda$-subharmonic function on
$T$ implying that $T$ is stochastically incomplete. \qed
\newline

This result has the following corollary:

\begin{corollary}\label{comp1cor}
If $T$ is a tree with a vertex $x_0$ such that
$\underline{m}_{x_0}(r) = \underline{m}(r) = \min_{x \in S_r(x_0)}
m(x)$ satisfies $\underline{m}(r) > 1$ and
\[ \sum_{r=0}^\infty \frac{1}{\underline{m}(r)} < \infty \]
then $T$ is stochastically incomplete. \end{corollary}

\Proof From the assumption on $T$, we can embed $T_n \subseteq T$,
where $T_n$ is defined by
\[  n(r) = \underline{m}(r) - 1  \textrm{ for } r>0 \]
and $n(0) = m(x_0)$.  Then $\sum_{r=0}^\infty \frac{1}{n(r)} <
\infty$ giving that $T_n$ is stochastically incomplete and so is,
therefore, $T$. \qed

\begin{remark}
This theorem and its corollary are unsatisfactory in the sense that
they require the tree to grow very rapidly in all directions from $x_0$
in order to be stochastically incomplete.  However, as we will see
in Theorem \ref{gen.tree.theorem} in the next section, it is sufficient that
the tree grows very rapidly in some direction from $x_0$.
\end{remark}

We now prove the inverse of the last result for a general graph
$G$.

\begin{theorem} \label{stochastic completeness}
If $G$ is any graph with a vertex $x_0$ such that $M_{x_0}(r) =
M(r) = \max_{x \in S_r(x_0)} m(x)$ satisfies
\[ \sum_{r=0}^\infty \frac{1}{M(r)} = \infty \]
then $G$ is stochastically complete.
\end{theorem}
\Proof Let $u$ be a positive, $\lambda$-harmonic function on $G$ for
$\lambda < 0$.  We will show that under the assumption on $G$, $u$
must be unbounded. At $x_0$, the relation $\Delta u(x_0) = \lambda
u(x_0)$, gives that
\begin{equation}\label{lamhar}
\sum_{x \sim x_0} u(x) = \big( m(x_0) - \lambda \big) u(x_0).
\end{equation}
This implies that there exists $x_1 \sim x_0$ such that
\[ u(x_1) \geq \left( 1 - \frac{\lambda}{m(x_0)} \right) u(x_0). \]
If not, then for all $x \sim x_0,$ $u(x) < \left( 1 -
\frac{\lambda}{m(x_0)} \right) u(x_0),$ giving that
\[ \sum_{x \sim x_0} u(x) < m(x_0)\left( 1 - \frac{\lambda}{m(x_0)}
\right) u(x_0) \] contradicting (\ref{lamhar}).

Now, by repeating the argument at $x_1$, we get that there must exist
a neighbor $y \sim x_1$ such that
\[ u(y) \geq \left( 1 - \frac{\lambda}{m(x_1)} \right) u(x_1). \]
Although $y$ is not necessarily in $S_2(x_0)$ we can repeat the argument
until we obtain a vertex $x_2 \in S_2(x_0)$ such that
\begin{eqnarray*}
u(x_2) &\geq& \left( 1 - \frac{\lambda}{m(x_1)} \right) u(x_1) \\
&\geq& \left( 1 - \frac{\lambda}{m(x_1)} \right) \left( 1 -
\frac{\lambda}{m(x_0)} \right) u(x_0).
\end{eqnarray*}

Iterating this argument, we get a sequence of distinct vertices $x_0
\sim x_1 \sim x_2 \sim \dots$ such that $x_r \in S_r(x_0)$ and
\[ u(x_r) \geq \prod_{i=0}^{r-1} \left( 1 - \frac{\lambda}{m(x_i)}
\right) u(x_0). \] Since $\sum_{i=0}^\infty \frac{1}{m(x_i)} \geq
\sum_{i=0}^{\infty}\frac{1}{M(i)} =\infty$ implies that
$\prod_{i=0}^\infty \left(1 - \frac{\lambda}{m(x_i)}\right) =
\infty$ it follows that $u$ cannot remain bounded, giving that $G$
must be stochastically complete. \qed

\begin{remark}
Theorem \ref{stochastic completeness} is a significant improvement over the result mentioned in the introduction \cite[Theorem 2.10]{Do-Mat} which states that the same conclusion as above holds if the valence is bounded above by a constant.  In fact the proof there can be extended to show that if $M(r)$ is $o(r)$ then the graph is stochastically complete whereas our result says that $M(r)$ can even be $O(r)$.
\end{remark}

A corollary of Theorem \ref{stochastic completeness} for trees is the following:
\begin{corollary}
Assume that $T \subseteq T_n$ with $x_0 \in T$.  If $T_n$ is
stochastically complete then so is $T$.
\end{corollary}
\Proof  Since $T_n$ is stochastically complete we have that
$\sum_{r=0}^\infty \frac{1}{n(r)} = \infty$ implying
$\sum_{r=0}^\infty \frac{1}{M(r)} = \infty$ so that $T$ is
stochastically complete.  \qed


\section{General Trees}
The purpose of this section is to follow-up on the remark following
Corollary \ref{comp1cor}.  The result there states that a general
tree $T$ will be stochastically incomplete if, starting out at a
fixed vertex, the branching grows rapidly in all directions.
The next theorem states that the same conclusion holds if the
branching grows rapidly in just one direction.

We start by slightly altering the notation used in the previous section.
If $x_0$ and $x_1$ are vertices of $T$ with $x_0 \sim x_1$ then we now denote
\[ \underline{m}(r) = \underline{m}_{\{x_0,x_1\}}(r) =
\min_{\substack{x \in S_r(x_0) \\ d(x,x_1) = r-1}} m(x) \quad \textrm{ for } r \geq 1\]
so that the minimum is now taken over those $x$ in $S_r(x_0)$ such
that $d(x,x_1) = r-1$.

\begin{theorem}\label{gen.tree.theorem}
If $T$ is a tree with a vertex $x_0 \in T$ such
that for some $x_1 \sim x_0$, $\underline{m}(r) =
\underline{m}_{\{x_0,x_1\}}(r)$ satisfies  $\underline{m}(r) > 1$ and
\[ \sum_{r=1}^{\infty} \frac{1}{\underline{m}(r)} < \infty \]
then $T$ is stochastically incomplete.
\end{theorem}

The proof of the theorem will use the following general proposition.
\begin{proposition} \label{gen.tree.prop}
For a graph $G$ with $x_0 \in G$ and any $\lambda < 0$ there exists
a function $v$ on $G$ such that $v(x_0)=1$, $0 < v(x) \leq 1$ for
all vertices $x$ and $\Delta v(x) = \lambda v(x)$ for all vertices
$x \not = x_0$.
\end{proposition}

\noindent \textbf{Proof of Theorem \ref{gen.tree.theorem}:} Assuming
Proposition \ref{gen.tree.prop}, we apply it to define a
positive, bounded $\lambda$-harmonic function $v$
on the part of $T$ below $x_0$ in Figure 3.1 with $v(x_0)=1$.
\begin{figure} \label{fig1} \centering
\includegraphics{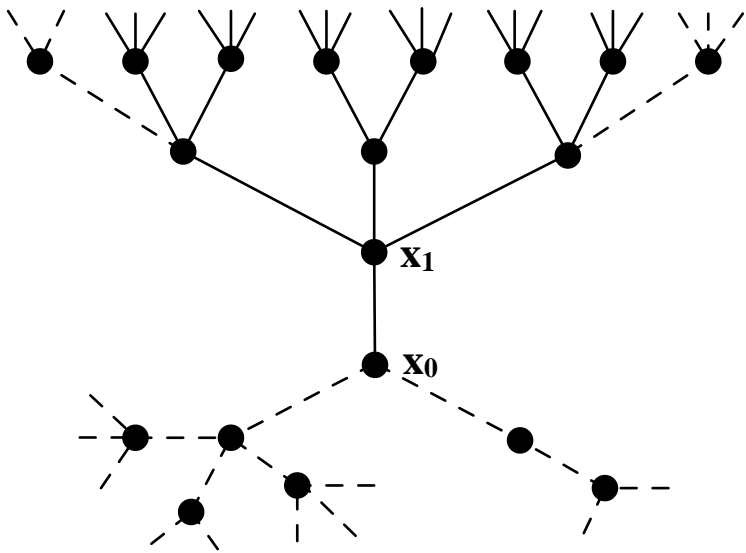}
\caption{$T$ with $T_n$ in solid.}
\end{figure}  For the part of $T$
above $x_0$ in Figure 3.1, the assumption on $T$ implies that we can
embed a stochastically incomplete model subtree $T_n$ with root
vertex $x_0$ where $n(0)=1$ and $n(r) =\underline{m}(r)-1$ for $r \geq 1$.
We extend the function $v$ to be defined on $T_n$ by
first making it $\lambda$-harmonic at $x_0$:
\[ v(x_1) = \big( m(x_0) - \lambda \big) - \sum_{\substack{x \sim
x_0 \\ x \not= x_1}} v(x). \]  Then, if $r$ denotes the distance to
the root, $v$ can be defined on the rest of $T_n$ by
\[ v(r+1) = \left(1 + \frac{1-\lambda}{n(r)}\right) v(r) -
\frac{1}{n(r)} v(r-1) \quad \textrm{ for } r \geq 2 \label{v(r)}.
\] This function will remain bounded since
$\sum_{r=1}^{\infty} \frac{1}{n(r)} < \infty.$ Therefore, by the
argument used in the proof of Theorem \ref{comp1}, there exists a
positive, bounded $\lambda$-subharmonic function on $T$. By
Theorem \ref{stochastic.incomp}, $T$ is then stochastically incomplete.  \qed \\

\noindent \textbf{Proof of Proposition \ref{gen.tree.prop}:} Let
$B_r(x_0)$ denote the ball of radius $r$ about $x_0$ in $G$.  For
$\lambda< 0$, on each $B_r(x_0)$ there exists a unique solution to
the following system of equations:
\begin{equation}\label{system of equations}
\left\{ \begin{array}{ll} \Delta v_r(x) = \lambda v_r(x) & \textrm{
for all }
x \in \textrm{int } B_r \setminus \{x_0\} \\
v_r(x_0)=1  \\
v_r(x) = 0 & \textrm{ for all } x \in \partial B_r.
\end{array} \right.
\end{equation}

Indeed, from basic linear algebra, since the system has the same
number of equations as there are values for $v_r$, there will exist
a unique solution if 0 is the only function which satisfies
\begin{equation}\label{homogeneous system}
\left\{  \begin{array}{ll}
\Delta v_r(x) =v_r(x) & \textrm{ for all } x \in \textrm{int } B_r \setminus \{x_0\}\\
v_r(x_0)=0  \\
v_r(x) = 0 & \textrm{ for all } x \in \partial B_r.
\end{array} \right.
\end{equation}

To show that this is so, suppose that $v_r$ is a non-zero solution
to (\ref{homogeneous system}). We can then assume that there
exists a vertex $\hat{x}$ in the interior such that $v_r(\hat{x}) >
0$ and $\hat{x}$ is a maximum for $v_r$ on $B_r(x_0)$. Then, by
calculation, $\Delta v_r(\hat{x}) \geq 0$, while $\Delta
v_r(\hat{x}) = \lambda v_r(\hat{x}) < 0$ giving a contradiction. The
same argument could be used to show that $v_r$ cannot have a
negative minimum.  Therefore, any solution to (\ref{homogeneous system}) must be zero.
This gives existence and uniqueness of a solution to (\ref{system of equations}).

Therefore, for $\lambda < 0$, for each $r$ there exists a unique
solution to
\[ \left\{  \begin{array}{ll}
\Delta v_r(x) = \lambda v_r(x) & \textrm{ for all } x \in \textrm{int }
 B_r \setminus \{x_0\} \\
v_r(x_0)=1  \\
v_r(x) = 0 & \textrm{ for all } x \in \partial B_r.
\end{array} \right. \]
On the interior of $B_r$, a solution must satisfy
\[0 < v_r \leq 1.\]
For, supposing that there exists a vertex $\hat{x}$ in the interior
such that $v_r(\hat{x}) \leq 0$ and $\hat{x}$ is a minimum for $v_r$
then, as before, by calculation, $\Delta v_r(\hat{x}) \leq 0$ while
$\Delta v_r(\hat{x}) = \lambda v_r(\hat{x}) \geq 0$. Therefore,
$v_r(x) = v_r(\hat{x})$ for all $x$ next to $\hat{x}$ and, by
repeating the argument, it would follow that $v_r$ is a constant
function for a non-positive constant yielding a contradiction since
$v_r(x_0) = 1$.  Hence, $v_r > 0$ implying $\Delta v_r < 0$ for
all vertices in the interior except for $x_0$ so that $v_r \leq 1$
since at an interior maximum $\Delta v_r \geq 0$.

Similarly, since $\Delta(v_{r+1} - v_r) = \lambda (v_{r+1}-v_r)$, it
follows that $v_r \leq v_{r+1}$ on  $B_r$ and by extending each
$v_r$ to be 0 outside of $B_r$, we get that
\[ v_r \leq v_{r+1} \textrm{ on } G.\]
Therefore, we can define $v$ as the limit
\[ v_r \to v \textrm{ as } r \to \infty.\]
It follows that $v$ satisfies $0 < v \leq 1$, $v(x_0) = 1$, and
$\Delta v = \lambda v$ for all vertices of $G$ except for $x_0$.
\qed


\section{Heat Kernel Comparison} The purpose of this section is to
prove two theorems which compare the heat kernel on a general tree
to the heat kernel on a model. These theorems were inspired by an
analogous result of Cheeger and Yau on model manifolds \cite[Theorem
3.1]{Ch-Y}.  Fixing a vertex $x_0$ in a tree $T$, we now denote
\[ \underline{m}(r) = \min_{x \in S_r(x_0)} m(x) \ \textrm{ and } \ M(r)
= \max_{x \in S_r(x_0)} m(x) \] the minimum and maximum valence
along the spheres $S_r(x_0)$.  Throughout, we use the notation
$\rho_t(x_0,x)$ for the heat kernel on $T_n$, while $p_t(x_0,x)$ will
denote the heat kernel on $T$. We will first show that, as a
function of $x$, $\rho_t(x_0,x)$ is constant on the spheres
$S_r(x_0)$ in $T_n$. Let $\rho_t(r) = \rho_t(0,r(x))$ denote this common
value.  Then the two main theorems of the section can be stated as
follows:

\begin{theorem}\label{heat.comp1}
If $M(r) \leq n(r) + 1$ for all $r>0$ then
\[ \rho_t(r) \leq p_t(x_0,x) \]
for all $x \in S_r(x_0) \subset T.$
\end{theorem}

\begin{theorem}\label{heat.comp2}
If $n(r) \leq \ \underline{m}(r)-1$ for all $r > 0$ then
\[p_t(x_0,x) \leq \rho_t(r) \]
for all $x \in S_r(x_0) \subset T$.
\end{theorem}

The proofs will follow easily from the maximum principle for the
heat equation once we establish two general lemmas
concerning the heat kernel on $T_n$. We start by proving the
property of the heat kernel mentioned at the start of this section.

\begin{lemma}
On $T_n$
\[ \rho_t(x_0,x) = \rho_t(r) \]
for all $x \in S_r(x_0).$
\end{lemma}

\Proof This result is essentially a restatement of the fact that the
coefficients of the Laplacian depend only on the valence which, on
$T_n$, only depends on the distance from the root.

We establish the result for the heat kernels $\rho_t^R(x_0,x)$ on
$B_R(x_0)$ with Dirichlet boundary conditions and pass to the limit.
The heat kernel $\rho_t^R(x_0,x)$ is the kernel of the operator
semigroup

\[ e^{-t \Delta_R} = I - t \Delta_R + \frac{t^2 \Delta_R^2}{2} -
\frac{t^3 \Delta_R^3}{6} + \dots \] where $\Delta_R$ denotes the
reduced Laplacian on $B_R(x_0)$.  That is,
\[ \rho_t^R(x_0,x) = \la \delta_{x_0},\delta_x \ra - t \Delta_R(x_0,x) +
\frac{t^2}{2} \Delta_R^2(x_0,x) - \dots \] where the coefficients of
the Laplacian $\Delta_R(x_0,x)$ are given by
\begin{eqnarray*}
\Delta_R(x_0,x) &=&  \Delta_R \delta_{x_0}(x) \\
&=& \left\{  \begin{array}{lll}
 n(0) & \textrm{ if } x = x_0 \\
 - 1 & \textrm{ if } x \in S_1(x_0) \\
 0  & \textrm{ otherwise }
 \end{array} \right.
\end{eqnarray*} and

\[ \Delta_R^{m+n}(x_0,x) = \sum_{y \in B_R} \Delta_R^m(x_0,y) \Delta_R^n(y,
x). \]  Therefore, these only depend on the distance between $x_0$ and
$x$. \qed

\begin{lemma}\label{rho2}
On $T_n$
\[ \rho_t(r) \geq \rho_t(r+1) \quad \textrm{ for all } r \geq 0. \]
\end{lemma}

\Proof We start with the general fact that for any graph
\[ \pdt p_t(x,x) \leq 0 \quad \textrm{ for all } t > 0.\]
Indeed, working with the heat kernels $p_t^R(x,y)$ on $B_R(x_0)$
with Dirichlet boundary conditions but now using the
eigenfunction expansion we have
\[ p_t^R(x,x) = \sum_{j=0}^{k(R)} e^{-\lambda_j^R t}
\big( \phi_j^R(x) \big)^2 \] implying
\[ \pdt p_t^R(x,x) = \sum_{j=0}^{k(R)} -\lambda_j^R
e^{-\lambda_j^R t} \big( \phi_j^R(x) \big)^2 < 0 \] since
$\lambda_j^R > 0$.  By passing to the limit we get that $\pdt
p_t(x,x) \leq 0.$

Therefore, in particular, $\pdt \rho_t^R(0) < 0$ which implies that
\[ \Delta \rho_t^R(0) = n(0) \Big( \rho_t^R(0) - \rho_t^R(1) \Big) > 0. \]
Thus,
\[ \rho_t^R(0) > \rho_t^R(1). \]
We also have that
\[ \rho_t^R(R-1) > \rho_t^R(R) = 0 \]
so that, as a function of $r$, $\rho_t^R(r)$ is decreasing at $r=0$
as well as at $r=R-1$ for all $t>0$.

Fix now a time $t_0>0$ and assume that there exists an $i_0>0$ such
that
\[ \rho_{t_0}^R(i_0) < \rho_{t_0}^R(i_0+1) \] where $i_0$ is the
smallest number with this property.  Therefore, as a function of
$r$, $\rho_{t_0}^R(r)$ achieves a local minimum at $r=i_0$. Since
$\rho_{t_0}^R(r)$ is decreasing again at $r=R-1,$ it follows that there
must exist a $j_0,$ $i_0<j_0 \leq R-1$ such that
\[ \rho_{t_0}^R(j_0) > \rho_{t_0}^R(j_0 + 1). \]
Therefore, $\rho_{t_0}^R(r)$ has a local maximum at $r=j_0$.  It
follows from calculation that
\[ \Delta \rho_{t_0}^R(i_0) < 0 \quad \textrm{ and } \quad \Delta
\rho_{t_0}^R(j_0) > 0 \] implying
\[ \pdt \rho_{t_0}^R(i_0) > 0 \quad \textrm{ and } \quad \pdt
\rho_{t_0}^R(j_0) < 0. \]
Therefore, at a previous time
$t_1 < t_0,$ $\rho_{t_1}^R(r)$, as a function of $r$, achieves a
smaller minimum than $\rho_{t_0}^R(i_0)$ and a larger maximum then
$\rho_{t_0}^R(j_0)$. That is, there exist $i_1$ and $j_1$ such that
\[ \rho_{t_1}^R(i_1) < \rho_{t_0}^R(i_0) \textrm { and }
\rho_{t_1}^R(j_1) > \rho_{t_0}^R(j_0) \] implying, in particular, that
\[ \rho_{t_1}^R(i_1) < \rho_{t_1}^R(j_1). \]

Since this argument can be repeated for any positive time $t$, it
follows that eventually we reach $t=0$ at which point we would have
an $i$ and $j$, $0 < i < j < R,$ such that
\[ \rho_0^R(i) < \rho_0^R(j) \]
contradicting the fact that $p_0^R(r) = 0$ for all $r \not = 0$.
Therefore, $\rho_t^R(r) \geq \rho_t^R(r+1)$ for all $r$ and
letting $R \to \infty$ we get
\[\rho_t(r) \geq \rho_t(r+1). \] \qed

\noindent \textbf{Proof of Theorems \ref{heat.comp1} and
\ref{heat.comp2}:}$\quad$   The proofs of the the two Theorems are
nearly identical so we give details for the proof of Theorem
\ref{heat.comp1} and then point out the modifications needed for the
proof of Theorem \ref{heat.comp2}.  We will denote the Laplacian on
$T_n$ by $\Delta_{T_n}$ to distinguish it from $\Delta_T$, the
Laplacian on $T$, wherever it is necessary.

For Theorem \ref{heat.comp1} it follows from the assumption $M(r)
\leq n(r)+1$ that we can embed $T$ into $T_n$.  Hence, we think of
$T \subseteq T_n$, and want to show that $\rho_t(r(x)) \leq
p_t(x_0,x)$ for all vertices of $T$, where $r(x) = d(x,x_0)$. We
work with the Dirichlet heat kernels on $B_R(x_0) \subset T$ and
consider the function
\[ u^R(x,t) = \rho_t^R(r(x)) - p_t^R(x_0,x)  \] on $B_R \times [0,S]$.
Then
\[ u^R(x,0) = 0 \textrm{ for all } x \in B_R \]
and \[ u^R(x,t) = 0 \textrm{ for all } x \in \partial B_R, \textrm{
all } t.
\]

Furthermore, it follows from $m(x) - 1 \leq n(r(x))$ and from Lemma
\ref{rho2} that $\rho_t^R(r(x))$ satisfies the following inequality
for any $x \in B_R$:
\begin{eqnarray*}
\Delta_T \rho_t^R(r(x)) &=& (m(x)-1)\Big(\rho_t^R(r(x)) -
\rho_t^R(r(x)+1)\Big) \\
& & + \ \rho_t^R(r(x)) - \rho_t^R(r(x)-1) \\
&\leq& n(r(x))\Big(\rho_t^R(r(x)) - \rho_t^R(r(x)+1)\Big)\\
& & + \ \rho_t^R(r(x)) - \rho_t^R(r(x) -1) \\
&=& \Delta_{T_n}\rho_t^R(r(x)) = -\pdt \rho_t^R(r(x)).
\end{eqnarray*}
Therefore, $u^R(x,t)$ satisfies
\[ \Delta_T u^R(x,t) + \pdt u^R(x,t) \leq 0.
\]
By applying the maximum principle for the
heat equation (see Remark \ref{maxheatremark}), it follows that
\[ \max_{B_R \times [0,S]} u^R(x,t) = \max_{\substack{B_R \times
\{0\} \ \cup \\ \partial B_R \times [0,S]}} u^R(x,t) = 0. \]
Therefore, $\rho_t^R(r(x)) - p_t^R(x,x_0) \leq 0$ so that  $\rho_t^R(r(x)) \leq p_t^R(x,x_0)$.
Since this holds for every $R$, by letting $R \to \infty$, we get that
\[ \rho_t(r(x)) \leq p_t(x,x_0) .\]

For Theorem \ref{heat.comp2}, from the assumption that $n(r) \leq \
\underline{m}(r) -1$, we may assume that $T_n \subseteq T$ and we
want to show that $p_t(x_0,x) \leq \rho_t(r(x))$. First, extend
$\rho$ to be defined on all of $T$ as before by letting:
\[ \rho_t(x) = \rho_t(r(x)) \textrm{ for } x \in T. \]
Then, since $n(r) \leq m(x) - 1$ for all $x \in
S_r \subset T$, it follows that
\[ u^R(x,t) = \rho_t^R(r(x)) - p_t^R(x_0,x) \] now satisfies
\[ \Delta_T u^R(x,t) + \pdt u_t^R(x,t) \geq 0. \]
This implies that
\[ \min_{B_R \times [0,S]} u^R(x,t) = \min_{\substack{B_R \times
\{0\} \ \cup \\ \partial B_R \times [0,S]}} u^R(x,t) = 0 \] implying
\[ p_t^R(x_0,x) \leq \rho_t^R(r(x)). \] \qed \\

In fact, using the same proof as above, we can extend the result of
Theorem \ref{heat.comp2} to a slightly more general graph $G$ in which
a model subtree $T_n$ may be embedded.  Essentially, we obtain $G$ by allowing any two
vertices on the same sphere to be connected by an edge in $T$ from Theorem \ref{heat.comp2}.
To make this precise, we introduce the following notation which will
be useful later as well.  Let $x_0 \in G$ be a fixed
vertex and $r(x) = d(x,x_0)$.
\begin{definition}
For $x \in G$ let
\begin{eqnarray*}
m_0(x) &=& \# \{y \ | \ y \sim x \textrm{ and } r(y) = r(x)\} \\
m_{+1}(x) &=& \# \{y \ | \ y \sim x \textrm{ and }r(y) = r(x)+1\} \\
m_{-1}(x) &=& \# \{y \ | \ y \sim x \textrm{ and }r(y) = r(x)-1\}
\end{eqnarray*}
\end{definition}  That is, $m_0(x),m_{+1}(x)$, and $m_{-1}(x)$ denote the number of
vertices that are the same distance, further away, and closer to
$x_0$ than is $x$, respectively.  We then state and prove the
following:
\begin{theorem}
If $G$ is any graph with $n(r) \leq m_{+1}(x)$ for all $x \in
S_r(x_0)$ and $m_{-1}(x) = 1$ for all vertices of $G$ then
\[ p_t(x_0,x) \leq \rho_t(r) \]
for all $x \in S_r(x_0) \subset G$.
\end{theorem}
\Proof  As before, we may assume that $T_n \subset G$.  Extend
$\rho$ to be defined on $G$ by letting:
\[ \rho_t(x) = \rho_t(r(x)) \textrm{ for } x \in G. \]]  Then, for $x \in S_r(x_0)
\subseteq G$,
\begin{eqnarray*}
\Delta_G \rho_t(r) &=& m_0(x)\big(\rho_t(r) - \rho_t(r) \big)
+ \ m_{+1}(x) \big( \rho_t(r) - \rho_t(r + 1) \big) \\
& & + \ m_{-1}(x) \big( \rho_t(r) - \rho_t(r - 1) \big) \\
&=& m_{+1}(x) \big( \rho_t(r) - \rho_t(r + 1) \big) + \ \rho_t(r) - \rho_t(r - 1) \\
&\geq& \Delta_{T_n} \rho_t(r) = -\pdt \rho_t(r).
\end{eqnarray*}
The rest of the proof is identical to the proof of Theorem
\ref{heat.comp2}. \qed


\section{Bounded Laplacian}
In this section, we introduce the bounded Laplacian $\Delta_{bd}$
and prove that, with this operator, any graph is stochastically
complete.  That is, in particular, bounded solutions to the heat
equation involving $\Delta_{bd}$ with bounded initial conditions are
unique. We refer to \cite{Do-Karp, Do-Ken} for the definitions
involved.

We define the bounded Laplacian to be the operator
\begin{eqnarray*}
\Delta_{bd}f(x) &=& f(x) - \frac{1}{m(x)}\sum_{y \sim x}f(y) \\
&=& \frac{1}{m(x)} \Delta f(x).
\end{eqnarray*}
In order for the analogue of Green's Theorem to hold we alter the inner
product on the space of functions on the graph.  We now let
\[ \la f,g \ra_{bd} = \sum_{x \in V} f(x) g(x) m(x) \]
while keeping the inner product on edges the same as before.  It is
now true that
\[ \la \Delta_{bd}f, g \ra_{bd} = \la df, dg \ra \]
for all $f$ such that
\[ \sum_{x \in V} f(x)^2 m(x) < \infty. \]

What distinguishes $\Delta_{bd}$ from $\Delta$ is that it is a bounded operator
without the assumption $m(x) \leq M$ necessary to imply that
$\Delta$ is bounded. This can be seen as follows.
\begin{eqnarray*}
\la df, df \ra &=& \sum_{[x,y] \in \tilde{E}}\big(f(y) -  f(x)\big)^2 \\
&\leq& 2 \sum_{[x,y] \in \tilde{E}}\big(f^2(y) + f^2(x) \big)\\
&=& 2 \sum_{x \in V}f^2(x)m(x) \ = \ 2 \la f,f \ra_{bd}.
\end{eqnarray*}
Therefore, $\lVert d \lVert \leq \sqrt 2$ implying $\lVert \Delta_{bd} \lVert \leq 2$.

We next prove that any graph is stochastically complete with respect
to this operator by studying $\lambda$-harmonic functions of
$\Delta_{bd}$.

\begin{theorem} \label{bounded.lap.thm}
If $v$ is a positive function on $G$ satisfying
\[\Delta_{bd}v(x) = \lambda v(x)\]
for $\lambda < 0$ then $v$ will be unbounded.
\end{theorem}

In particular, since the proof of the equivalence of the various
formulations of stochastic incompleteness for $\Delta$ holds for
$\Delta_{bd}$, we have that

\begin{corollary}
For any bounded function $u_0$, the bounded solution to
\[ \left\{  \begin{array}{ll} \Delta_{bd} u(x,t) + \pdtu(x,t) = 0 &
\textrm{ for all } x \in V, \textrm{all }  t>0\\
u(x,0) = u_0(x) & \textrm{ for all } x \in V
 \end{array} \right. \]
is unique.
\end{corollary}

\noindent \textbf{Proof of Theorem \ref{bounded.lap.thm}:}  The
proof is essentially the proof of Theorem \ref{stochastic
completeness} rewritten for the bounded Laplacian. Fix a vertex
$x_0$ of $G$.  We will show that there exists a sequence of distinct
vertices
\[ x_0 \sim x_1 \sim x_2 \sim \dots \]
such that
\[ v(x_i) \to \infty \textrm{ as } i \to \infty.\]
At $x_0$,
\[ \Delta_{bd} v(x_0) = v(x_0) - \frac{1}{m(x_0)} \sum_{x \sim x_0}
v(x) = \lambda v(x_0) \] implies that
\begin{equation}\label{bounded equality}
\sum_{x \sim x_0} v(x) = m(x_0)(1 - \lambda)v(x_0).
\end{equation}
Therefore, there exists a neighbor $x_1$ of $x_0$ such that
\[ v(x_1) \geq (1-\lambda) v(x_0).\]
Since, if not, if $v(x) < (1-\lambda)v(x_0)$ for all $x \sim x_0$, then
\[ \sum_{x \sim x_0} v(x) < m(x_0)(1-\lambda)v(x_0) \] contradicting
(\ref{bounded equality}).  Applying the argument now at $x_1$ we get a
neighbor $x_2$ of $x_1$ such that
\[ v(x_2) \geq (1-\lambda)v(x_1) \geq (1-\lambda)^2 v(x_0). \]
In general, we get a sequence of distinct vertices $x_0 \sim
x_1 \sim x_2 \sim \dots$ such that
\[ v(x_i) \geq (1-\lambda)^i v(x_0) \]
implying that
\[ v(x_i) \to \infty \textrm{ as } i \to \infty.\] \qed


\chapter{Spectral Analysis}

\section{Bottom of the Spectrum}
We recall the definition of $\lambda_0(\Delta)$, the bottom of the spectrum of
the Laplacian on a general graph $G$, and prove a characterization of
it in terms of $\lambda$-harmonic functions.  This result was
inspired by an analogous result in \cite[Theorem 2.1]{Sull} and was
proven for the bounded Laplacian in \cite{Do-Karp}.

Fix a vertex $x_0$ in $G$ and let $B_r= B_r(x_0)$ denote the ball of
radius $r$ about $x_0$ with boundary $\partial B_r$ as before. Also,
let $\Delta_r$ denote the reduced Laplacian acting on the space
$C(B_r, \partial B_r)$ of functions on $B_r$ that vanish on the
boundary $\partial B_r$.  We define then $\lambda_0^r =
\lambda_0(\Delta_r)$ as
\[ \lambda_0^r = \lambda_0(\Delta_r) = \min_{\substack{f \in C(B_r, \partial
B_r) \\ f \not \equiv 0}} \frac{\la df, df \ra}{\la f,f \ra}
\]
and show, as in \cite[Lemma 1.9]{Do}, that
\begin{lemma}\label{Ray-lem}
$\lambda_0^r$ is the smallest eigenvalue of $\Delta_r$.
Furthermore, if $f_0$ is a function in $C(B_r,
\partial B_r)$ such that
\begin{equation}\label{Ray}
\lambda_0^r = \frac{\la df_0, df_0 \ra}{\la f_0,f_0 \ra}
\end{equation}
then $\Delta_r f_0 = \lambda_0^r f_0$ and $f_0$ can be chosen so
that $f_0 > 0$ on the interior of $B_r$.
\end{lemma}
\Proof  If $\lambda$ is any eigenvalue of $\Delta_r$ with
eigenfunction $f$ then
\[ \frac{\la df,df \ra}{\la f,f \ra} = \frac{\la \Delta_r f,f \ra}{\la f,f \ra} = \lambda \]
implies that $\lambda \geq \lambda_0^r$.

Now, if $f_0$ satisfies (\ref{Ray}) above and $\{\lambda_i^r
\}_{i=0}^{k(r)}$ are the eigenvalues of $\Delta_r$ with
$\{\phi_i^r\}_{i=0}^{k(r)}$ a set of corresponding eigenfunctions which are
an orthonormal basis for $C(B_r,\partial B_r)$ then
\[ f_0 = \sum_{i=0}^{k(r)} a_i \phi_i^r \]
where $a_i = \la f_0, \phi_i^r \ra$.  We wish to show that $a_i=0$
if $\lambda_i^r \not = \lambda_0^r$.  This can be seen as follows
{\setlength\arraycolsep{0pt}
\begin{eqnarray*}
0 \ &\leq& \ \left\la d\Big(f_0 - \sum_{i=0}^{k(r)} a_i
\phi_i^r\Big),
d\Big(f_0 - \sum_{j=0}^{k(r)} a_j \phi_j^r\Big)  \right\ra\\
&=&\la df_0, df_0 \ra - 2 \sum_{i=0}^{k(r)} a_i \la f_0,\Delta_r
\phi_i^r \ra + \sum_{i,j = 0}^{k(r)} a_i a_j \la \phi_i^r, \Delta_r \phi_j^r \ra \\
&=& \la df_0,df_0 \ra - 2 \sum_{i=0}^{k(r)} a_i^2 \lambda_i^r +
\sum_{i,j=0}^{k(r)} a_i a_j \lambda_j^r \la \phi_i^r, \phi_j^r \ra  \\
&=& \la df_0,df_0 \ra - \sum_{i=0}^{k(r)} a_i^2 \lambda_i^r
\end{eqnarray*}}
implies that
\[ \la df_0,df_0 \ra \ \geq \ \sum_{i=0}^{k(r)} a_i^2 \lambda_i^r. \]
While (\ref{Ray}) gives
\[ \la df_0,df_0 \ra = \lambda_0^r \la f_0,f_0 \ra = \lambda_0^r
\sum_{i=0}^{k(r)} a_i^2.\]
Therefore, $a_i = 0$ if
$\lambda_i^r \not = \lambda_0^r$.

Now, noting that
\[ \la f_0, f_0 \ra = \la |f_0|,|f_0| \ra \]
while
\[ \la df_0, df_0\ra \ \geq \ \la d|f_0|,d|f_0| \ra \]
it is clear that (\ref{Ray}) can only be decreased by replacing
$f_0$ by $|f_0|$ and we may assume at the onset that $f_0 \geq 0$.
Then, if there exists a vertex $\hat x$ in the interior of $B_r$
where $f_0(\hat x) = 0$ then it follows from $\Delta_r f_0 =
\lambda_0^r f_0$ that
\[ \Delta_r f_0(\hat x) = -\sum_{x \sim \hat x}f_0(x) = 0.\]
Therefore, $f_0(x) = 0$ for all $x \sim \hat x$.  Repeating this
argument would give that $f_0 = 0$ on the interior of $B_r$ yielding
a contradiction. \qed \newline

It follows from Lemma \ref{Ray-lem} that
\[ \lambda_0^r \geq \lambda_0^{r+1} > 0 \]
so we may define
\[ \lambda_0 = \lambda_0(\Delta) = \lim_{r \to \infty} \lambda_0(\Delta_r). \]

\begin{remark}
It is clear that this number is independent of the choice  of
exhaustion sequence for the graph $G$ since if
$\{D_i\}_{i=0}^{\infty}$ is any other exhaustion sequence then for
each $R$ there is $I_R$ large enough so that $B_R \subset D_{I_R}$.
Therefore, by Lemma \ref{Ray-lem}, $\lambda_0^R \geq \lambda_0^{D_{I_R}}$.
Reversing the roles of $B_r$ and $D_i$, we get that $\lambda_0^{B_r}$ and $\lambda_0^{D_i}$
converge to the same number.  Also, for future reference, we point out that
$\lambda_0(\Delta)$ can also be defined as
\[ \lambda_0(\Delta) = \inf_{\substack{f \in C_0(V) \\ f \not \equiv 0}} \frac{\la
df, df \ra}{\la f, f \ra} \] where $C_0(V)$ denotes the set of finitely
supported functions on the graph $G$.
\end{remark}

We now state and prove the following characterization of
$\lambda_0(\Delta)$ in terms of $\lambda$-harmonic functions
\cite[Theorem 2.1]{Sull}.

\begin{theorem} \label{lambda0}
For every $\lambda \leq \lambda_0(\Delta)$ there exists a positive
$\lambda$-harmonic function.  For every $\lambda >
\lambda_0(\Delta)$ there is no such function.
\end{theorem}

\Proof  The proof of the first part is a variation of the
argument given in \cite[Theorem 2.4]{Do}.  See also
\cite[Proposition 1.5]{Do-Karp} for the case of the bounded
Laplacian and \cite[Lemma 1]{Fis-Sh} for manifolds.

We start with the case of $\lambda = \lambda_0(\Delta).$ From Lemma
\ref{Ray-lem}, for each $r$ there exists a positive function $v_r$
such that
\[ \Delta_r v_r = \lambda_0^r v_r \textrm{ on int } B_r. \]
We normalize this function so that $v_r(x_0) = 1$ and extend it to
be 0 outside of $B_r$. We will show that this function is bounded
for all vertices $x$. Let $M(i) = M_{x_0}(i) = \max_{x \in S_i(x_0)}
m(x)$ as before.  Then if $x_1 \in S_1 \subset$ int $B_r$ it follows
from $\Delta v_r(x_0) > 0$ that
\[ m(x_0) v_r(x_0) > \sum_{x \sim x_0} v_r(x) \geq v_r(x_1) \]
implying
\[ v_r(x_1) < M(0). \]
By repeating the same argument we get that if $x_i \in S_i$ where $i
< r$ then
\begin{equation} \label{upperbound}
v_r(x_i) < M(i-1) M(i-2) \ldots M(0).
\end{equation}

Now, using the diagonal process we can find a subsequence of $\{ v_r
\}_{r=1}^{\infty}$ which converges for all vertices $x$.  Denote this as
\[ v_{r_k}(x) \to v(x) \textrm{ as } k \to \infty \textrm{ for all } x.\]
It follows that $\Delta v(x) = \lambda_0 v(x)$ for all vertices $x$,
$v \geq 0$ and $v(x_0) = 1.$  By using the same argument as in the proof of Lemma
\ref{Ray-lem} if there exists an $x$ where $v(x) = 0$ then $v$ would
have to be constantly 0 yielding a contradiction since $v(x_0)=1$.
This completes the case of $\lambda = \lambda_0.$

For the case of $\lambda < \lambda_0(\Delta)$ we modify the argument
above as follows.  First, as noted in the proof of the implication $2') \Rightarrow 2)$
in Theorem \ref{stochastic.incomp}, since $\lambda < \lambda_0(\Delta) \leq
\lambda_0(\Delta_r)$ the operator $(\Delta_r - \lambda I)$ is
positive and hence invertible on $C(B_r, \partial B_r)$, the space of all functions on $B_r$ which vanish on $\partial B_r$, so that one can find a function $v_r$
which satisfies
\[ \left\{  \begin{array}{ll} \Delta_r v_r = \lambda v_r & \textrm{on int } B_r\\
 {v_r}_{| \partial B_r} = 1.
 \end{array} \right. \]
Indeed, as before, one can let $v_r = (\Delta_r - \lambda
I)^{-1}(\lambda_{\textrm{int }B_r}) + \bf{1}$, where $\lambda_{\textrm{int }B_r}$ is the
function equal to $\lambda$ on every vertex in the interior of $B_r$
and 0 elsewhere and \bf{1} is equal to 1 on every vertex of $B_r$.
We renormalize $v_r$ so that it is equal to 1 at $x_0$ and call it
$u_r$, that is, let
\[u_r = \frac{1}{v_r(x_0)} v_r.\]

Now, if $u_r \geq 0$ on the interior of $B_r$ then $u_r > 0$ on the
interior of $B_r$ by the same argument as above. To show that $u_r
\geq 0$ we assume that there exists a vertex $\hat x$ in the
interior of $B_r$ where $u_r(\hat x) < 0$ and let $w$ be a new
function defined by
\[w(x) = \left\{  \begin{array}{ll} u_r(x) & \textrm{for $x$ such that } u_r(x) < 0\\
   0  & \textrm{otherwise}.
 \end{array} \right. \]
In particular, $w(x) = 0$ if $x$ is in the boundary of $B_r$, so the
error term vanishes when we apply the analogue of Green's Theorem to
$w$ below.  Now, if $x$ is a vertex in the interior of $B_r$ such
that $u_r<0$ for $x$ and all neighbors of $x$ then $\Delta w(x) =
\Delta u_r(x)$. If $x$ is a vertex where $u_r(x) < 0$ and $x$ has a
neighbor $y$ for which $u_r(y) \geq 0$ then $\Delta w(x) \geq \Delta
u_r(x)$.  Combining these, it follows that
\begin{eqnarray*}
\la dw,dw \ra &=& \la \Delta w, w \ra \\
&=& \sum_{\substack{x \in \textrm{int } B_r \\ u_r(x) < 0}} \Delta w(x) w(x) \\
&\leq& \sum_{\substack{x \in \textrm{int } B_r \\ u_r(x) < 0}} \Delta u_r(x) u_r(x) \\
&=& \lambda \sum_{\substack{x \in \textrm{int } B_r \\ u_r(x) < 0}} u_r(x)^2 \\
&=& \lambda \la w, w \ra
\end{eqnarray*}
so that
\[\frac{\la dw, dw \ra}{\la w, w \ra} \leq \lambda.\]
From Lemma \ref{Ray-lem} it would follow that $\lambda_0(\Delta_r)
\leq \lambda$ contradicting the assumption that $\lambda <
\lambda_0(\Delta)$. Therefore, $u_r \geq 0$ on the interior of $B_r$
and so $u_r > 0$ there as well.

Now, if $0 < \lambda < \lambda_0(\Delta)$ then we can use the same
argument as above to show that $u_r(x)$ is bounded for all vertices
$x$ as in (\ref{upperbound}). If $\lambda \leq
0$ then from $\Delta u_r = \lambda u_r$ the bound becomes
\[ u_r(x_1) \leq M(0) - \lambda \]
for all $x_1 \in S_1$ and
\[ u_r(x_i) \leq (M(i-1) - \lambda) (M(i-2) - \lambda) \dots (M(0) - \lambda)\]
for all $x_i \in S_i$.  The remainder of the argument is the same as
before. This completes the proof of the first part of Theorem
\ref{lambda0}.

The proof of the second part of Theorem \ref{lambda0} is adapted
from \cite[$\textrm{p.}$ 761]{Sull}.  See \cite{Do-Karp} for a different proof
involving the use of Green's Theorem.  Suppose that there exists
a positive function $v$ such that $\Delta v = \lambda v$. Then,
letting
\[ u(x,t) = e^{-\lambda t}v(x) \]
and
\[ w(x,t) = \sum_{y \in B_r}p_t^r(x,y) v(y)\]
we see that both $u$ and $w$ satisfy the heat
equation on int $B_r \times (0,S)$ and $u(x,0)=w(x,0) = v(x)$ on int
$B_r$ with $u(x,0) = v(x) > 0$ and $w(x,0) = 0$ on $\partial B_r$.
By applying the maximum principle for the heat equation to the
difference of the two functions we get that
\[ \min_{B_r \times [0,S]}(u-w) = \min_{\substack{B_r \times \{0\}
\cup\\ \partial B_r \times [0,S]}} (u-w) \geq 0\] by what was noted
above and since $w(x,t)$ vanishes on $\partial B_r$ while $u(x,t)$
is positive there.  Therefore, $u(x,t) \geq w(x,t)$ or
\[e^{-\lambda t}v(x) \geq \sum_{y \in B_r} p_t^r(x,y)v(y) \textrm{ on } B_r \times [0,S] .\]
By using the eigenvalue and eigenfunction expansion for $p_t^r(x,y)$
we now get that
\begin{eqnarray*}
v(x) &\geq& \sum_{y \in B_r}e^{\lambda t} p_t^r(x,y) v(y) \\
&=& \sum_{y \in B_r}e^{\lambda t} \sum_{i=0}^{k(r)}e^{-\lambda_i^r
t} \phi_i^r(x) \phi_i^r(y) v(y) \\
&=& \sum_{y \in B_r}\sum_{i=0}^{k(r)} e^{(\lambda - \lambda_i^r)t}
\phi_i^r(x) \phi_i^r(y) v(y)
\end{eqnarray*}
and if $\lambda > \lambda_0^r$ then the right hand side would tend
to $\infty$ as $t \to \infty$.  It follows that $\lambda \leq
\lambda_0^r$ for all $r$ so that $\lambda \leq \lambda_0(\Delta)$.
\qed


\subsection{Relationship to Stochastic Incompleteness}
In light of the large volume growth that is required for a graph to
be stochastically incomplete and the well-known relationship between
the bottom of the spectrum and Cheeger's constant which is, at least
partially, outlined in this subsection and the next section, it
might seem plausible to conjecture that stochastic incompleteness
would imply that $\lambda_0(\Delta) > 0$.  The purpose of this
subsection is to give an example where this is not the case.

The example is constructed as follows: start with a model tree $T_n$
which is stochastically incomplete, that is, such that
$\sum_{r=0}^{\infty} \frac{1}{n(r)} < \infty$, and attach to the
root vertex $x_0$ an infinitely long path $x_0 \sim x_1 \sim x_2
\sim \ldots.$  The resulting tree will be stochastically incomplete
by Theorem \ref{gen.tree.theorem}.

We now show that $\lambda_0(\Delta) = 0$.  As noted before
\begin{equation}\label{lambda0bd}
\lambda_0 \leq \frac{\la df, df \ra}{\la f, f \ra}
\end{equation}
for any nonzero, finitely supported function $f$.  By taking any
finite subgraph $D$ and substituting its characteristic function
$1_D$ into (\ref{lambda0bd}) we get that
\[ \lambda_0 \leq \frac{L(\partial D)}{\textrm{Vol} (D)} \]
where
\[ L(\partial D) = \# \{y \sim x \  | \ x \in D \textrm{ and } y \not \in D
\} \] that is, the number of edges with one vertex in $D$ and one
not in $D$, and Vol$(D)$ denotes the number of vertices in $D$.  By
taking increasingly larger connected subgraphs of the path that was
added onto out model tree $T_n$, it is clear that this ratio goes to
0 since $L(\partial D) = 2$, for all such subgraphs.


\section{Lower Bounds}
In this section we prove some estimates of $\lambda_0(\Delta),$ the
bottom of the spectrum of $\Delta.$ In order to prove our results,
we work with the bounded Laplacian $\Delta_{bd}$ and use the
characterization of $\lambda_0(\Delta_{bd})$ in terms of Cheeger's
constant proved in \cite{Do-Ken} to get a lower bound for $\lambda_0(\Delta_{bd})$.
We then transfer this result to obtain a lower bound for $\lambda_0(\Delta)$.

We recall that $\Delta_{bd}$ is given by
\[ \Delta_{bd}f(x) = f(x) - \frac{1}{m(x)}\sum_{y \sim x}f(y) =
\frac{1}{m(x)} \Delta f(x).
\]
The bottom of the spectrum is, as for $\Delta$, given by an exhaustion
argument or, equivalently, as
\[ \lambda_0(\Delta_{bd}) = \inf_{\substack{ f \in C_0(V) \\ f \not \equiv 0 }}
\frac{\la \Delta_{bd} f,f\ra_{bd}}{\la f,f \ra_{bd}} \] where the
infimum is taken over all nonzero, finitely supported functions $f$ and the inner
product is now given by
\[ \la f, f \ra_{bd} = \sum_{x \in V}f(x)^2 m(x). \]

For a finite subgraph $D$ we define $A(D),$ the area of $D$, to be
\[ A(D) = \sum_{x \in D} m(x) \] and $L(\partial D)$, the length of the boundary, to
be
\[ L(\partial D) = \# \{y \sim x \ | \ x \in D \textrm{ and } y \not \in D
\} \] as in the previous subsection.  If we let
\[ \alpha = \inf_{\substack{D \subset G \\ D \textrm{ finite, connected }}} \frac{L(\partial
D)}{A(D)} \] then the main result in Section 2 of \cite{Do-Ken}
states that
\begin{theorem}\label{Do-K}
\[ \lambda_0(\Delta_{bd}) \geq \frac{\alpha^2}{2}.\]
\end{theorem}

Now, for a general graph $G$ we fix a vertex $x_0$ of $G$ and let
$r(x) = d(x,x_0)$.  Then for a vertex $x$ we let
\begin{eqnarray*}
m_0(x) &=& \# \{y \ | \ y \sim x \textrm{ and } r(y) = r(x) \} \\
m_{+1}(x) &=& \# \{y \ | \ y \sim x \textrm{ and } r(y) = r(x) + 1\} \\
m_{-1}(x) &=& \# \{y \ | \ y \sim x \textrm{ and } r(y) = r(x) - 1\}
\end{eqnarray*}
as before.  The main theorem of this section can now be stated as
follows:
\begin{theorem}\label{lam0lbthm}
If for all vertices $x$ of $G$
\[ \frac{m_{+1}(x) - m_{-1}(x) }{m(x)} \geq c >0   \]
then
\[ \lambda_0(\Delta_{bd}) \geq \frac{c^2}{2}. \]
If, in addition, $m(x) \geq m$ then
\begin{equation} \label{lam0lb}
\lambda_0(\Delta) \geq \frac{c^2}{2} m.
\end{equation}
\end{theorem}

\begin{example}
\textup{For a tree we have that $m_0(x) = 0$ and $m_{-1}(x) = 1$ for
all vertices $x$ so that}
\[ \frac{m_{+1}(x) - m_{-1}(x)}{m(x)} = \frac{m(x)-2}{m(x)} = 1 - \frac{2}{m(x)}.
\]
\textup{Therefore, if $m(x) \geq m > 2$ then $\frac{m_{+1}(x) -
m_{-1}(x) }{m(x)} \geq \left( 1 - \frac{2}{m} \right)$. Hence,
for such a tree, we get $c = \frac{m-2}{m}$, implying}
\[ \lambda_0(\Delta_{bd}) \geq \frac{(m-2)^2}{2m^2} \]
\textup{and}
\[ \lambda_0(\Delta) \geq \frac{(m-2)^2}{2m}. \]
\end{example}

\Proof  Take $D \subset G$ a finite, connected subgraph.  Let $r(x) = d(x,x_0)$
where $x_0$ is a fixed vertex of $G$. Then by applying the analogue
of Green's Theorem we get
\begin{eqnarray}\label{L ineq}
\left| \sum_{x \in D} \Delta_{bd}r(x)m(x) \right| &=& \left| \sum_{x
\in D} \Delta r(x) \right|  \nonumber \\
&=& \left| \sum_{\substack{y \sim x \\ x \in D, y \not \in D}}
\big(r(x) - r(y) \big) \right|  \nonumber \\
&\leq& \sum_{\substack{y \sim x \\ x \in D, y \not \in D}}
|r(x) - r(y)| \nonumber \\
&\leq& L(\partial D)
\end{eqnarray}
since $r(x) - r(y)$ can only be $\pm 1$ or $0$ if $y \sim x$.

On the other hand,
\begin{eqnarray*}
\Delta_{bd}r(x) &=&  r(x) - \frac{1}{m(x)}\Big( m_0(x)r(x) +
m_{+1}(x)\big(r(x) + 1\big) + m_{-1}(x)\big(r(x) - 1\big) \Big) \\
\\
&=&  \ \frac{m_{-1}(x) - m_{+1}(x)}{m(x)}
\end{eqnarray*}
since $m(x) = m_0(x) + m_{+1}(x) + m_{-1}(x)$.  Therefore, it
follows from the assumption $ \frac{m_{+1}(x) - m_{-1}(x) }{m(x)} \geq c $
that $\Delta_{bd}r(x) < 0 $ and
\[ \left| \Delta_{bd}r(x) \right| \geq c \]
for all vertices $x$.  Therefore,
\begin{eqnarray}\label{A ineq}
\left| \sum_{x \in D} \Delta_{bd}r(x)m(x) \right| &=&  \sum_{x \in
D} | \Delta_{bd}r(x)m(x) |  \nonumber \\
&\geq& c \sum_{x \in D} m(x) = c A(D).
\end{eqnarray}
Combining the inequalities (\ref{L ineq}) and (\ref{A ineq}), we get
\[ cA(D) \leq L(\partial D) \]
or that
\[ c \leq \frac{L(\partial D)}{A(D)} \]
for all finite, connected subgraphs $D$.  Applying Theorem
\ref{Do-K} it follows that
\begin{equation}\label{lam0bd}
 \lambda_0(\Delta_{bd}) \geq \frac{c^2}{2}.
\end{equation}  This gives the
first part of Theorem \ref{lam0lbthm}.

For the second part of Theorem \ref{lam0lbthm}, we proceed as
follows. By using the Rayleigh-Ritz characterization of
$\lambda_0(\Delta_{bd})$, inequality (\ref{lam0bd}) gives
\[ \la \Delta_{bd}f,f \ra_{bd} \ \geq \frac{c^2}{2} \la f,f \ra_{bd} \]
for every finitely supported, nonzero function on the graph $G$.
Now,
\begin{eqnarray*}
\la \Delta_{bd}f,f \ra_{bd} &=& \sum_{x \in V} \Delta_{bd}f(x) f(x)
m(x)
\\
&=& \sum_{x \in V} \Delta f(x) f(x) \ = \ \la \Delta f, f \ra
\end{eqnarray*}
while, if $m(x) \geq m$, then
\begin{eqnarray*}
\frac{c^2}{2} \la f,f \ra_{bd} &=& \frac{c^2}{2} \sum_{x \in V} f(x)^2 m(x) \\
&\geq& \frac{c^2}{2} m \sum_{x \in V} f(x)^2 \ = \ \frac{c^2}{2} m
\la f,f \ra.
\end{eqnarray*}
Therefore, we get
\[ \frac{\la \Delta f, f \ra}{\la f,f \ra} \geq \frac{c^2}{2} m \]
and taking the infimum over the set of all finitely supported,
nonzero functions $f$ it follows that
\[ \lambda_0(\Delta) \geq \frac{c^2}{2} m.  \] \qed


\section{Essential Spectrum}
We now use Theorem \ref{lam0lbthm} to prove that, under certain
assumptions on the graph, $\tilde{\Delta}$, the unique self-adjoint extension
of $\Delta$ to $\ell^2(V)$, has empty
essential spectrum as in \cite[Theorem 1.1]{Don-Li}.

The essential spectrum is, by definition, the complement in the
spectrum of the set of isolated eigenvalues of finite multiplicity.
We use the notation spec($\tilde{\Delta}$) and ess spec($\tilde{\Delta}$) for the
spectrum and essential spectrum of $\tilde{\Delta}$ respectively.  Now, as
pointed out in \cite[Theorem VII.12 and remarks following Theorem
VIII.6]{Reed-Simon}, the essential spectrum of a self-adjoint
operator can be characterized as follows
\begin{theorem}
$\lambda \in$  \textup{ess spec}$(\tilde{\Delta})$ if and only if there
exists an sequence of orthonormal function $\{f_i\}_{i=0}^{\infty}$
in the domain of $\tilde{\Delta}$ such that
\[ \lim_{i \to \infty} \| \tilde{\Delta} f_i - \lambda f_i \|_{\ell^2} = 0. \]
\end{theorem}

In fact, it is sufficient that the sequence be noncompact, that is, have
no convergent subsequence, and this will be the characterization of
the essential spectrum that we use to prove the theorem below.

Fix a vertex $x_0$ and let $\underline{m}(r)$ denote the
smallest valence of the vertices on the sphere $S_r(x_0)$ as before:
\[ \underline{m}(r) = \min_{x \in S_r(x_0)} m(x). \]  We then have

\begin{theorem} \label{essspecthm}
If for all vertices $x$ of $G$
\[ \frac{m_{+1}(x) - m_{-1}(x) }{m(x)} \geq c > 0\]
and
\[ \underline{m}(r) \to \infty \textrm{ as }  r \to \infty\] then
$\tilde{\Delta}$ has empty essential spectrum.
\end{theorem}

\begin{example}
Again, for a tree, we note that if $m(x) \geq m > 2$ for all
vertices $x$ then the first assumption is satisfied and so if
$\underline{m}(r) \to \infty$ then $\tilde{\Delta}$ will have empty
essential spectrum.
\end{example}

The proof of Theorem \ref{essspecthm} will follow easily once we
establish the following lemma which is analogous to
\cite[Proposition 2.1]{Don-Li} and apply the second result of
Theorem \ref{lam0lbthm}.  Let $\tilde{\Delta}_r$ denote the self-adjoint
extension of the Laplacian acting on the space $C_0(V,B_r)$, of
functions with finite support disjoint from $B_r$, to
$\ell^2(V,B_r)$, the square summable functions which vanish on
$B_r$. We then have that
\begin{lemma}\label{essspeclem}
$\tilde{\Delta}$ and $\tilde{\Delta}_r$ have the same essential spectrum.
\end{lemma}

Assuming the lemma for now, we give the proof of the theorem:

\noindent \textbf{Proof of Theorem \ref{essspecthm}: }  By applying
(\ref{lam0lb}) from Theorem \ref{lam0lbthm} we get that
\[ \lambda_0(\tilde{\Delta}_r) \to \infty \textrm{ as } r \to \infty \]
since $\underline{m}(r) \to \infty$.  Now, applying Lemma
\ref{essspeclem}, since the essential spectrum of $\tilde{\Delta}$ is the
same as that of $\tilde{\Delta}_r$ and the bottom of the spectrum of
$\tilde{\Delta}_r$ is increasing to infinity, it must follow that the
essential spectrum of $\tilde{\Delta}$ is empty. \qed \newline

\noindent \textbf{Proof of Lemma \ref{essspeclem}: } Let $\lambda
\in $ ess spec $(\tilde{\Delta}_r)$.  Let $\{f_i\}_{i=0}^\infty$ be a
sequence of orthonormal functions vanishing on $B_r$ satisfying
\[ \lim_{i \to \infty} \| \tilde{\Delta}_r f_i - \lambda f_i \|_{\ell^2} = 0.  \]
Then, since
\[ \tilde{\Delta}_r f_i(x) \not = \tilde{\Delta} f_i(x) \textrm{ only for } x \in
\partial B_r \]
and, by orthonormality, for every vertex $x$, $f_i(x) \to 0$ as $i
\to \infty$ it follows that
\[ \lim_{i \to \infty} \| \tilde{\Delta} f_i - \lambda f_i \|_{\ell^2} = 0  \]
so that $\lambda \in $ ess spec $(\tilde{\Delta})$.

Now, say that $\lambda \in \textrm{ess spec}(\tilde{\Delta})$ and
$\{f_i\}_{i=0}^\infty$ is a sequence of orthonormal function such
that
\[ \| \tilde{\Delta} f_i - \lambda f_i \|_{\ell^2} \to 0 \textrm{ as } i \to \infty.
\]
Let
\[ \varphi_r(x)= \left\{ \begin{array}{ll} 0 & \textrm{if } x \in
B_r(x_0) \\
1 & \textrm{otherwise} \end{array} \right.
\]
We claim that $\{ \varphi_r f_i \}_{i=0}^{\infty}$ will be a
sequence of bounded functions with no convergent subsequence
satisfying
\[ \| \tilde{\Delta}_r (\varphi_r f_i) - \lambda (\varphi_r f_i) \|_{\ell_2} \to 0
\textrm{ as } i \to \infty. \]

To show that $\{ \varphi_r f_i \}_{i=0}^\infty$ has no convergent
subsequences we first note that since $\{ f_i \}_{i=0}^\infty$ are
orthonormal,  $\{ f_i \}_{i=0}$ has no convergent subsequences.  This
follows since pointwise, using orthonormality as above,
\[ f_i(x) \to 0 \textrm{ as } i \to \infty \ \textrm{ for all } x \in V\]
while $\| f_i \|_{\ell^2} = 1 $ for all $i$.  Now, assume that $\{
\varphi_r f_i \}_{i=0}^\infty$ has a convergent subsequence, say
\[ \varphi_r f_{i_k} \to f \textrm{ as } k \to \infty \textrm{ in } \ell_2. \]
Since $f_i \in \ell_2$, $\{ f_{i_k}(x) \}_{k=0}^\infty$ has a
convergent subsequence for each $x$.  Because $B_r$ has only finitely
many vertices, we can find a subsequence of $\{ f_{i_k}
\}_{k=0}^\infty$ which converges for each $x \in B_r$.  We continue
to denote this subsequence as $\{ f_{i_k} \}_{k=0}^\infty$ and let
$\hat{f}(x)$ be defined by
\[ \hat{f}(x) =  \left\{ \begin{array}{ll} f(x) & \textrm{ if } x \not \in B_r \\
\lim_{k \to \infty}  f_{i_k}(x)  & \textrm{ if }x \in B_r
\end{array} \right.
\]
Then, it would follow that
\begin{eqnarray*}
\| f_{i_k} - \hat{f} \|^2_{\ell_2} &=& \sum_{x \in V}\big(
f_{i_k}(x) - \hat{f}(x) \big)^2 \\
&=& \sum_{x \not \in B_r}\big(f_{i_k}(x) - f(x) \big)^2 + \sum_{x
\in B_r}\big(f_{i_k}(x) - \hat{f}(x)\big)^2 \\
&\to& 0 \textrm{ as } k \to \infty.
\end{eqnarray*}
so that $\{ f_i \}_{i=0}^\infty$ would have a convergent
subsequence.  The contradiction shows that $\{ \varphi_r f_i
\}_{i=0}^\infty$ cannot have a convergent subsequence.

What remains to be shown is that $\| (\tilde{\Delta}_r - \lambda I)
(\varphi_r f_i) \|_{\ell^2} \to 0$ as $i \to \infty$.  First, we
calculate $\tilde{\Delta}_r (\varphi_r f_i)$:
\begin{eqnarray*}
\tilde{\Delta}_r (\varphi_r f_i)(x) &=& \sum_{y \sim x}\big(
\varphi_r(x)f_i(x) - \varphi_r(y)f_i(y) \big) \\
&=& \sum_{y \sim x}\big( \varphi_r(x)f_i(x) - \varphi_r(x)f_i(y) +
\varphi_r(x)f_i(y) - \varphi_r(y)f_i(y) \big) \\
&=& \varphi_r(x) \tilde{\Delta}_r f_i(x) + \sum_{y \sim x}f_i(y)
\big(\varphi_r(x) - \varphi_r(y) \big).
\end{eqnarray*}
Therefore, {\setlength\arraycolsep{0pt}
\begin{eqnarray*}
\sum_{x \in V}\big(\tilde{\Delta}_r(\varphi_r f_i)(x) - \lambda (\varphi_r
f_i)(x) \big)^2 &\leq& \sum_{x \in V} \Big(\varphi_r(x) \big(
(\tilde{\Delta}_r f_i)(x) - \lambda f_i(x)\big) \Big)^2 \\
&\quad & + \sum_{x \in V}\sum_{y \sim x} \Big(f_i(y)
\big(\varphi_r(x) -
\varphi_r(y) \big) \Big)^2 \\
&=& \sum_{x \not \in B_r} \big( (\tilde{\Delta} f_i)(x) -
\lambda f_i(x)\big)^2 \\
&\quad& + \sum_{x \in \partial B_r}\sum_{\substack{y \sim x \\ y
\not \in B_r}}\big( f_i(x)^2 + f_i(y)^2 \big).
\end{eqnarray*}}

Now, the first sum above goes to 0 as $i \to \infty$ by the
assumption on $f_i$ while, since $f_i$ are orthonormal, $f_i(x) \to
0$ as $i \to \infty$ for each $x$, so it follows that the second sum
also goes to 0.  Therefore,
\[ \| \tilde{\Delta}_r (\varphi_r f_i) - \lambda (\varphi_r f_i) \|_{\ell_2} \to 0
\textrm{ as } i \to \infty \] and so $\lambda \in $ ess
spec$(\tilde{\Delta}_r)$. \qed


\end{document}